# DISCRETE APPROXIMATIONS TO REFLECTED BROWNIAN MOTION[1]

By Krzysztof Burdzy and Zhen-Qing Chen

*University of Washington*

In this paper we investigate three discrete or semi-discrete approximation schemes for reflected Brownian motion on bounded Euclidean domains. For a class of bounded domains $D$ in $\mathbb{R}^n$ that includes all bounded Lipschitz domains and the von Koch snowflake domain, we show that the laws of both discrete and continuous time simple random walks on $D \cap 2^{-k}\mathbb{Z}^n$ moving at the rate $2^{-2k}$ with stationary initial distribution converge weakly in the space $\mathbf{D}([0,1], \mathbb{R}^n)$, equipped with the Skorokhod topology, to the law of the stationary reflected Brownian motion on $D$. We further show that the following "myopic conditioning" algorithm generates, in the limit, a reflected Brownian motion on any bounded domain $D$. For every integer $k \geq 1$, let $\{X^k_{j2^{-k}}, j = 0, 1, 2, \ldots\}$ be a discrete time Markov chain with one-step transition probabilities being the same as those for the Brownian motion in $D$ conditioned not to exit $D$ before time $2^{-k}$. We prove that the laws of $X^k$ converge to that of the reflected Brownian motion on $D$. These approximation schemes give not only new ways of constructing reflected Brownian motion but also implementable algorithms to simulate reflected Brownian motion.

**1. Introduction.** Let $n \geq 1$ and $D \subset \mathbb{R}^n$ be a domain (connected open set) with compact closure. Consider a reflected Brownian motion (RBM in abbreviation) $Y$ in $D$. Heuristically, RBM in $D$ is a continuous Markov process $Y$ taking values in $\overline{D}$ that behaves like a Brownian motion in $\mathbb{R}^n$ when $Y_t \in D$ and is instantaneously pushed back along the inward normal direction when $Y_t \in \partial D$. RBM on smooth domains can be constructed in various ways, including the deterministic Skorokhod problem method, martingale problem method, or as a solution to a stochastic differential equation with reflecting boundary conditions (see the Introduction of [5]). When $D$ is

Received November 2006; revised February 2007.
[1]Supported in part by NSF Grant DMS-06-0206.

nonsmooth, all the methods mentioned above cease to work. On nonsmooth domains, the most effective way to construct RBM is to use the Dirichlet form method. The RBM constructed using a Dirichlet form coincides with RBM constructed using any other standard method in every smooth domain.

It is natural to try to construct RBM in a nonsmooth domain using a sequence of approximations that can be easily constructed themselves. One such approximating scheme was studied in [2] and [5], where the processes approximating RBM in a nonsmooth domain $D$ were RBMs in smooth domains increasing to $D$. In this paper we will consider processes approximating RBM in $D$ that are defined on the same state space $D$, or a discrete subspace of $D$. We consider the subject interesting and important in itself, but we also have a more concrete and direct motivation—it comes from two recent papers on multi-particle systems [3, 4]. In each of these papers a large population of particles is trapped in a domain. All particles perform independent random walks or Brownian motions and an appropriate mechanism keeps all particles inside the domain, which is assumed to be regular in a certain sense (see [3, 4] for more details). The intuitive interpretation of those results and their extension to less regular domains would require an invariance principle for the reflected random walk and similar processes in domains with nonsmooth boundaries. As far as we can tell, such results are not available in literature and they do not follow easily from published theorems.

In this paper we investigate three discrete or semi-discrete approximation schemes for reflected Brownian motion. The first two approximations involve random walks and we prove that they converge to reflected Brownian motion in a class of bounded nonsmooth domains in $\mathbb{R}^n$ that includes all bounded Lipschitz domains and the von Koch snowflake domain. The third scheme is based on "myopic" conditioning and it converges to the reflected Brownian motion in all bounded domains. We will now describe these schemes in more detail.

Let $D$ be a bounded domain in $\mathbb{R}^n$ whose boundary $\partial D$ has zero Lebesgue measure. Without loss of generality, we may assume that $0 \in D$. Let $D_k$ be the connected component of $D \cap 2^{-k}\mathbb{Z}^n$ that contains 0 with edge structure inherited from $2^{-k}\mathbb{Z}^n$ (see the next section for a precise definition). We will use $v_k(x)$ to denote the degree of a vertex $x$ in $D_k$. Let $X^k$ and $Y^k$ be the discrete and continuous time simple random walks on $D_k$ moving at the rate $2^{-2k}$ with stationary initial distribution $m_k$, respectively, where $m_k(x) = \frac{v_k(x)}{2n} 2^{-kn}$. We show that the laws of both $\{X^k, k \geq 1\}$ and $\{Y^k, k \geq 1\}$ are tight in the Skorokhod space $\mathbf{D}([0,\infty), \mathbb{R}^n)$ of right continuous functions having left limits. We show (see Theorems 2.4 and 3.3 below) that if $D$ satisfies an additional condition (1.1) below, which is satisfied by all bounded Lipschitz domains and all bounded uniform domains (see below for the definition), then both $\{X^k, k \geq 1\}$ and $\{Y^k, k \geq 1\}$ converge weakly to the stationary reflected Brownian motion on $D$ in the Skorokhod space $\mathbf{D}([0,1], \mathbb{R}^n)$.

The last of our main theorems is concerned with "myopic conditioning." We say that a Markov process is conditioned in a myopic way if it is conditioned not to hit the boundary for a very short period of time, say, $2^{-k}$ units of time, where $k$ is large. At the end of this period of time, we restart the process at its current position and condition it to avoid the boundary for another period of $2^{-k}$ units of time. We repeat the conditioning step over and over again. Intuition suggests that, when $2^{-k}$ is very small and the process is far from the boundary, the effect of conditioning is negligible. On the other hand, one expects that when the process is very close to the boundary, the effect of conditioning is a strong repulsion from the boundary. These two heuristic remarks suggest that, for small $2^{-k}$, the effect of myopic conditioning is similar to that of reflection. A more precise description of myopic conditioning of Brownian motion is the following. For every integer $k \geq 1$, let $\{Z^k_{j2^{-k}}, j = 0, 1, 2, \ldots\}$ be a discrete time Markov chain with one-step transition probabilities being the same as those for the Brownian motion in $D$ conditioned not to exit $D$ before time $2^{-k}$. The process $Z^k_t$ can be defined for $t \in [(j-1)2^{-k}, j2^{-k}]$ either as the conditional Brownian motion going from $Z^k_{(j-1)2^{-k}}$ to $Z^k_{j2^{-k}}$ without leaving the domain $D$ or as a linear interpolation between $Z^k_{(j-1)2^{-k}}$ and $Z^k_{j2^{-k}}$. We prove in Theorems 5.1 and 5.6 below that, for any bounded domain $D$, the laws of $Z^k$ (defined in either way) converge to that of the reflected Brownian motion on $D$. We remark that, in Theorems 5.1 and 5.6, the myopic conditioning approximation of reflected Brownian motion is proved for every starting point $x \in D$, so these theorems demonstrate explicitly that the symmetric reflected Brownian motion on $D$ is completely determined by the absorbing Brownian motion in $D$.

We would like to point out that for the first two approximation schemes (i.e., random walk approximations) discussed in this paper, the proof of tightness of the approximating sequences requires only the weak assumption that $D$ is bounded and its boundary has zero Lebesgue measure. However, an example given in Section 4 shows that in some domains $D$, reflected random walks do not converge to the reflected Brownian motion in $D$.

The definition of the random walk on a lattice depends very much on the geometric structure of the state space. Myopic conditioning, on the other hand, can be applied to any Markov process that has a positive probability of not hitting a set for any fixed amount of time, if it starts outside that set. Hence, myopic conditioning might provide a new way of defining reflected Markov processes, for example, the reflected stable process that was introduced in [1], whenever the limit can be shown to exist. We plan to address this problem in a future work.

The literature on topics discussed in this paper is rather limited. An invariance principle for discrete approximations to the reflected Brownian motion was proved in [15] in $C^2$-domains. Our approach in Theorem 2.4 of this paper is quite different from that in [15]. For results on approximations to the *killed* Brownian motion in Lipschitz domains, see, for example, [18].

In the rest of this introduction we give a brief review of RBM on nonsmooth domains, as well as present some definitions and notation, followed by a brief description of the approach of this paper to establish discrete approximations of reflected Brownian motion.

Let $n \geq 1$ and $D$ be a bounded connected open set in $\mathbb{R}^n$. Denote by $m$ the Lebesgue measure in $\mathbb{R}^n$. Define

$$W^{1,2}(D) := \{f \in L^2(D,m) : \nabla f \in L^2(D,m)\},$$

equipped with norm $\|f\|_{1,2} := \|f\|_2 + \|\nabla f\|_2$, where $\|f\|_p$ denotes the $L^p(D,m)$-norm of $f$ for $p \geq 1$. In Sections 2 and 3 we will assume that the boundary of $D$ has zero volume, that is, $m(\partial D) = 0$. To prove the weak convergence of random walks, we need to impose the following condition on $D$:

(1.1) $\qquad\qquad C^1(\overline{D})$ is dense in $(W^{1,2}(D), \|\cdot\|_{1,2})$.

Here $C^1(\overline{D})$ is the space of real-valued continuous functions on $\overline{D}$ that have continuous first derivatives on $\overline{D}$. Condition (1.1) is satisfied when $D$ is a $W^{1,2}$-extension domain in the sense that there is a bounded linear operator $T$ from $(W^{1,2}(D), \|\cdot\|_{1,2})$ to $(W^{1,2}(\mathbb{R}^n), \|\cdot\|_{1,2})$ so that $Tf = f$ $m$-a.e. on $D$ for every $f \in W^{1,2}(D)$. This is because $C_c^\infty(\mathbb{R}^n)$, the space of smooth functions with compact support in $\mathbb{R}^n$, is $\|\cdot\|_{1,2}$-dense in $W^{1,2}(\mathbb{R}^n)$. Examples of $W^{1,2}$-extension domains are bounded Lipschitz domains in $\mathbb{R}^n$, and, more generally, local uniform domains also known as $(\varepsilon, \delta)$-domains (see [12]), defined as follows. We say that $D$ is an $(\varepsilon, \delta)$-*domain* if $\delta, \varepsilon > 0$, and whenever $x, y \in D$ and $|x - y| < \delta$, then there exists a rectifiable arc $\gamma \subset D$ joining $x$ and $y$ with length$(\gamma) \leq \varepsilon^{-1}|x - y|$ such that $\min\{|x - z|, |z - y|\} \leq \varepsilon^{-1}\operatorname{dist}(z, \partial D)$ for all points $z \in \gamma$. Here $\operatorname{dist}(z, \partial D)$ is the Euclidean distance between a point $z$ and the set $\partial D$. "Uniform domains" can be defined as $(\varepsilon, \infty)$-domains. An example of a uniform domain is the classical von Koch snowflake domain. Every *nontangentially accessible domain* defined by Jerison and Kenig in [11] is a uniform domain (see (3.4) of [11]). The boundary of a uniform domain can be highly nonrectifiable and, in general, no regularity of its boundary can be inferred (besides the easy fact that the Hausdorff dimension of the boundary is strictly less than $n$). For any $\alpha \in [n-1, n)$, one can construct a uniform domain $D \subset \mathbb{R}^n$ such that $\mathcal{H}^\alpha(U \cap \partial D) > 0$ for any open set $U$ satisfying $U \cap \partial D \neq \varnothing$. Here $\mathcal{H}^\alpha$ denotes the $\alpha$-dimensional Hausdorff measure in $\mathbb{R}^n$.

Define for $f, g \in W^{1,2}(D)$ a bilinear form

$$\mathcal{E}(f,g) = \tfrac{1}{2} \int_D \nabla f(x) \cdot \nabla g(x) m(dx).$$

Under condition (1.1), the Dirichlet form $(\mathcal{E}, W^{1,2}(D))$ is regular and, therefore, there is a symmetric diffusion $X$ taking values in $\overline{D}$ associated with it, called the reflected Brownian motion on $D$. See [9] for definitions and properties of Dirichlet spaces, including the notions of quasi-everywhere, quasi-continuous, and so on. When $D$ is $C^1$-smooth, $X$ admits the following Skorokhod decomposition (cf. [5]):

$$X_t = X_0 + B_t + \int_0^t \mathbf{n}(X_s) \, dL_s, \qquad t \geq 0,$$

where $B$ is the standard Brownian motion in $\mathbb{R}^n$, $\mathbf{n}$ is the unit inward normal vector field on $\partial D$ and $L$ is a continuous nondecreasing additive functional of $X$ (called the boundary local time of $X$) that increases only when $X_t$ is on the boundary of $D$.

Constructing a reflected Brownian motion on a nonsmooth bounded domain $D$ is a delicate problem. Fukushima [8] used the Martin–Kuramochi compactification $D^*$ of $D$ to construct a continuous symmetric diffusion process $X^*$ on $D^*$ whose Dirichlet form is $(\mathcal{E}, W^{1,2}(D))$. The process $X^*$ could be called reflected Brownian motion in $D$, but it lives on an abstract space $D^*$ that contains $D$ as a dense open set. It was proposed in [5] to refer the quasi-continuous projection $X$ of $X^*$ from $D^*$ into the Euclidean closure $\overline{D}$ as reflected Brownian motion in $D$. The projection process $X$ is a continuous Markov process on $\overline{D}$ and it spends zero Lebesgue amount of time on $\partial D$. Moreover, it behaves like Brownian motion when it is inside $D$. For every point $x \in D$, one can construct reflected Brownian motion in $D$ defined as above so that it starts from $x$. The distribution of $X$ is uniquely determined by the fact that its associated Dirichlet form is $(\mathcal{E}, W^{1,2}(D))$. In general, $X$ is not a strong Markov process on $\overline{D}$ (e.g., this is the case when $D$ is the unit disk with a slit removed). However, when condition (1.1) is satisfied, $X$ is the usual reflected Brownian motion in $D$ obtained as the Hunt process associated with the regular Dirichlet form $(\mathcal{E}, W^{1,2}(D))$ on $\overline{D}$. See [2], Section 3 and [5], Section 1 for more information on the history of reflected Brownian motion on nonsmooth domains.

For any metric state space $\mathcal{S}$ and positive $T > 0$, let $\mathbf{D}([0,T], \mathcal{S})$ denote the space of all functions on $[0,T]$ taking values in $\mathcal{S}$ that are right continuous and have left limits. The space $\mathbf{D}([0,\infty), \mathcal{S})$ is a separable metric space under the Skorokhod topology. We refer the reader to [7] for the definition and properties of the Skorokhod topology and space $\mathbf{D}([0,T], \mathcal{S})$. The space of continuous functions from $[0,T]$ to $S$, equipped with the local uniform topology, will be denoted by $C([0,T], \mathcal{S})$. We will use similar notation for spaces of functions defined on finite time intervals. For $D \subset \mathbb{R}^n$, we will use $C_c^\infty(D)$ to denote the space of smooth functions with compact support in $D$.

For technical convenience, we will often consider stochastic processes whose initial distribution is a finite measure, not necessarily normalized to have total mass 1, for example, the Lebesgue measure on a bounded set $D$. Translating our results to the usual probabilistic setting is straightforward and so it is left to the reader.

We close this section with a brief description of the main ideas of our proofs. Denote by $\{P_t^D, t \geq 0\}$ the transition semigroup of Brownian motion killed upon leaving domain $D$. In the first two approximation schemes (i.e., random walk approximations), define $m_k(x) = \frac{v_k(x)}{2n} 2^{-kn}$ for $x \in D_k$, where $v_k(x)$ is the degree of the vertex $x$ in $D_k$; and in the myopic conditioning scheme, define $m_k(dx) := 1_D(x) P_{2^{-k}}^D 1(x) m(dx)$. Let $X^k$ be one of the discrete approximating processes mentioned above. Then $m_k$ is the reversible measure for $X^k$ in a suitable sense. We will show that the law of $\{X^k, \mathbf{P}_{m_k}^k, k \geq 1\}$ is tight in the space $C([0,1], \mathbb{R}^n)$ or $D([0,1], \mathbb{R}^n)$, and that any of its weak subsequential limits $(Z, \mathbf{P})$ is a time-homogeneous Markov process that is time-reversible with respect to the Lebesgue measure $m$ in $D$. We then show that the process $Z$ killed upon leaving domain $D$ is a killed Brownian motion in $D$ and establish that the Dirichlet form $(\mathcal{E}^Z, \mathcal{F})$ of $Z$ in $L^2(D, m)$ has the property that

$$W^{1,2}(D) \subset \mathcal{F}$$

and

$$\mathcal{E}^Z(f,f) \leq \mathcal{E}(f,f) := \tfrac{1}{2} \int_D |\nabla f(x)|^2 \, dx \qquad \text{for } f \in W^{1,2}(D).$$

We thus conclude from the following Theorem 1.1 that $(\mathcal{E}^Z, \mathcal{F}) = (\mathcal{E}, W^{1,2}(D))$ and, therefore, $(Z, \mathbf{P})$ is the stationary reflected Brownian motion on $D$. Discrete approximations with nonstationary initial distributions are then handled via those with stationary distributions.

THEOREM 1.1. *Let $D$ be a bounded domain in $\mathbb{R}^n$ and $m_D$ be the Lebesgue measure in $D$ that is extended to $\overline{D}$ by taking $m_D(\overline{D} \setminus D) = 0$. Suppose that $Z$ is a $\overline{D}$-valued right continuous time-homogeneous Markov process having left-limits with initial distribution $m_D$ and is symmetric with respect to measure $m_D$. Let $(\mathcal{E}^Z, \mathcal{F})$ be the Dirichlet form of $Z$ in $L^2(\overline{D}, m_D)$. If the subprocess of $Z$ killed upon leaving domain $D$ is a killed Brownian motion in $D$, then*

$$\mathcal{F} \subset W^{1,2}(D) \quad \text{and} \quad \mathcal{E}^Z(f,f) \geq \mathcal{E}(f,f) \qquad \text{for } f \in \mathcal{F}.$$

The above theorem is essentially due to Silverstein [13], Theorems 15.2 and 20.1. In fact, these theorems were proved in a more general context. Observe that the Dirichlet form for the killed Brownian motion in $D$ is $(\mathcal{E}, W_0^{1,2}(D))$ and $(\mathcal{E}, W^{1,2}(D))$ is its active reflected Dirichlet space. Here $W_0^{1,2}(D)$ is the completion of $C_c^\infty(D)$ under the norm $\|\cdot\|_{1,2}$. An accessible proof of Theorem 1.1 can be found in Takeda [17], Theorem 3.3.

**2. Invariance principle for discrete reflected random walk.** Let $D$ be a bounded connected open set in $\mathbb{R}^n$ with $m(\partial D) = 0$. Without loss of generality, assume that the origin $0 \in D$. Let $\overline{2^{-k}\mathbb{Z}^n}$ be the union of all closed line segments joining nearest neighbors in $2^{-k}\mathbb{Z}^n$, let $D_k^*$ be the connected component of $\overline{2^{-k}\mathbb{Z}^n} \cap D$ that contains the point 0, and let $D_k = D_k^* \cap 2^{-k}\mathbb{Z}^n$. For $x \in D_k$, we use $v_k(x)$ to denote the degree of the vertex $x$ in $D_k$. Let $\{X_{j2^{-2k}}^k, j = 0, 1, \ldots\}$ be the simple random walk on $D_k$ that jumps every $2^{-2k}$ unit of time. By definition, the random walk $\{X_{j2^{-2k}}^k, j = 0, 1, \ldots\}$ jumps to one of its nearest neighbors with equal probabilities. This discrete time Markov chain is symmetric with respect to measure $m_k$, where $m_k(x) = \frac{v_k(x)}{2n} 2^{-kn}$ for $x \in D_k$. Clearly, $m_k$ converge weakly to $m$ on $D$. We now extend the time-parameter of $\{X_{j2^{-2k}}^k, j = 0, 1, \ldots\}$ to all nonnegative reals using linear interpolation over the intervals $((j-1)2^{-2k}, j2^{-2k})$ for $j = 1, 2, \ldots$. We thus obtain a process $X^k = \{X_t^k, t \geq 0\}$. Its law with $X_0^k = x$ will be denoted by $\mathbf{P}_x^k$.

For $x, y \in D_k$, let $x \leftrightarrow y$ mean that $x$ and $y$ are at the distance $2^{-k}$. Let $Q_k(x, dy)$ denote the one-step transition probability for the discrete time Markov chain $\{X_{j2^{-2k}}^k, j = 0, 1, \ldots\}$; that is, for $f \geq 0$ on $D$ and $x \in D_k$,

$$Q_k f(x) := \int_D f(y) Q_k(x, dy) := \frac{1}{v_k(x)} \sum_{y \in D_k : y \leftrightarrow x} f(y).$$

For $f \in C^2(\overline{D})$, define

$$\mathcal{L}_k f(x) := \int_D (f(y) - f(x)) Q_k(x, dy)$$

$$= \frac{1}{v_k(x)} \sum_{y \in D_k : y \leftrightarrow x} (f(y) - f(x)), \qquad x \in D_k.$$

Then $\{f(X_{j2^{-2k}}^k) - \sum_{i=0}^{j-1} \mathcal{L}_k f(X_{i2^{-2k}}^k), \mathcal{G}_{j2^{-k}}^k, j = 0, 1, \ldots\}$ is a martingale for every $f \in C^2(\overline{D})$, where $\mathcal{G}_t^k := \sigma(X_s^k, s \leq t)$.

To study the weak limit of $\{X^k, k \geq 1\}$, we introduce an auxiliary process $Y^k$ defined by $Y_t^k := X_{[2^{2k}t]2^{-2k}}^k$, where $[\alpha]$ denotes the largest integer that is less than or equal to $\alpha$. Note that $Y^k$ is a time-inhomogeneous Markov process. For every fixed $t > 0$, its transition probability operator is symmetric with respect to the measure $m_k$ on $D_k$. Let $\mathcal{F}_t^k := \sigma(Y_s^k, s \leq t)$. By abuse of notation, the law of $Y^k$ starting from $x \in D_k$ will also be denoted by $\mathbf{P}_x^k$.

LEMMA 2.1. *Let $D$ be a bounded domain in $\mathbb{R}^n$ with $m(\partial D) = 0$. Then the laws $\{\mathbf{P}_{m_k}^k, k \geq 1\}$ of $\{X^k, k \geq 1\}$ are tight in the space $C([0,T], \mathbb{R}^n)$ for every $T > 0$.*

PROOF. For each fixed $k \geq 1$, we may assume, without loss of generality, that $\Omega$ is the canonical space $\mathbf{D}([0,\infty), \mathbb{R}^n)$ and $Y_t^k$ is the coordinate map on $\Omega$. Given $t > 0$ and a path $\omega \in \Omega$, the time reversal operator $r_t$ is defined by

$$(2.1) \qquad r_t(\omega)(s) := \begin{cases} \omega((t-s)-), & \text{if } 0 \leq s \leq t, \\ \omega(0), & \text{if } s \geq t. \end{cases}$$

Here for $r > 0$, $\omega(r-) := \lim_{s \uparrow r} \omega(s)$ is the left limit at $r$, and we use the convention that $\omega(0-) := \omega(0)$. We note that

$$(2.2) \qquad \begin{aligned} \lim_{s \downarrow 0} r_t(\omega)(s) &= \omega(t-) = r_t(\omega)(0) \quad \text{and} \\ \lim_{s \uparrow t} r_t(\omega)(s) &= \omega(0) = r_t(\omega)(t). \end{aligned}$$

Observe that for every integer $T \geq 1$, $\mathbf{P}_{m_k}^k$ restricted to the time interval $[0, T)$ is invariant under the time-reversal operator $r_T$. Note that

$$M_t^{k,f} := f(Y_t^k) - f(Y_0^k) - \sum_{i=0}^{[2^{2k}t]-1} \mathcal{L}_k f(Y_{i2^{-2k}}^k)$$

is an $\{\mathcal{F}_t^k, t \geq 0\}$-martingale for every $f \in C^2(\overline{D})$. We have

$$(2.3) \quad f(Y_t^k) - f(Y_0^k) = \tfrac{1}{2} M_t^{k,f} - \tfrac{1}{2}(M_{T-}^{k,f} - M_{(T-t)-}^{k,f}) \circ r_T \qquad \text{for } t \in [0,T).$$

For $f_i(x) = x_i$, let $M^{k,i} := M^{k,f_i}$ and $M^k := (M^{k,1}, \ldots, M^{k,n})$. With this notation, we have

$$Y_t^k - Y_0^k = \tfrac{1}{2} M_t^k - \tfrac{1}{2}(M_{T-}^k - M_{(T-t)-}^k) \circ r_T \qquad \text{for } t \in [0,T).$$

Let $S_k := \{x \in D_k : v(k) = 2n\}$. For $f \in C^3(\overline{D})$, by the Taylor expansion, for $x \in D_k$,

$$\mathcal{L}_k f(x) = \int_D \left( \sum_{i=1}^n \frac{\partial f(x)}{\partial x_i}(y_i - x_i) \right.$$
$$\left. + \frac{1}{2} \sum_{i,j=1}^n \frac{\partial^2 f(x)}{\partial x_i \partial x_j}(y_i - x_i)(y_j - x_j) + O(1)|y-x|^3 \right) Q_k(x, dy).$$

We see that, for $x \in S_k$,

$$\mathcal{L}_k f(x) = \int_D \left( \frac{1}{2} \sum_{i=1}^n \frac{\partial^2 f(x)}{\partial x_i^2} (y_i - x_i)^2 + O(1)|y-x|^3 \right) Q_k(x, dy),$$

and so

(2.4) $$\mathcal{L}_k f(x) = \frac{1}{2n} \Delta f(x) 2^{-2k} + O(1) 2^{-3k},$$

while for $x \in D_k \setminus S_k$,

$$|\mathcal{L}_k f(x)| \le 2^{-2k} O(1).$$

For $\delta > 0$, define $D^\delta := \{x \in \mathbb{R}^n : \text{dist}(x, \partial D) < 2\delta\}$. Note that $D^{2^{-k}\sqrt{n}}$ contains the $\sqrt{n} 2^{-k}$-neighborhood of $D_k \setminus S_k$ and so $m_k(D_k \setminus S_k) \le m(D^{2^{-k}\sqrt{n}})$, which goes to 0 as $k \to \infty$ because $m(\partial D) = 0$. Thus, for every $t > 0$,

$$[M^{k,i}, M^{k,j}]_t = \sum_{l=1}^{[2^{2k}t]} (Y^{k,i}_{l2^{-2k}} - Y^{k,i}_{(l-1)2^{-2k}} - \mathcal{L}_k x_i(Y^k_{(l-1)2^{-2k}}))$$

$$\times (Y^{k,j}_{l2^{-2k}} - Y^{k,j}_{(l-1)2^{-2k}} - \mathcal{L}_k x_j(Y^k_{(l-1)2^{-2k}}))$$

$$= \delta_{ij} \sum_{l=1}^{[2^{2k}t]} ((Y^{k,i}_{l2^{-2k}} - Y^{k,i}_{(l-1)2^{-2k}})^2$$

$$+ O(1) 2^{-3k} + O(1) 2^{-2k} \mathbf{1}_{\{Y^k_{(l-1)2^{-2k}} \in D_k \setminus S_k\}}).$$

Since $m_k$ is the invariant measure for the Markov chain $\{Y^k_{j2^{-2k}}, j = 0, 1, \ldots\}$, we have

$$\lim_{k\to\infty} \mathbf{E}^k_{m_k}\left[ \left| [M^{k,i}, M^{k,j}]_t - \delta_{ij} \frac{t}{n} \right| \right] \le \lim_{k\to\infty} tm(D^{2^{-k}\sqrt{n}}) O(1) = 0.$$

By [7], Theorem 7.4.1, the laws of $M^k$ converge weakly to that of a Brownian motion in $\mathbb{R}^n$ in the space $\mathbf{D}([0,T], \mathbb{R}^n)$. Since $\mathbf{P}^k_{m_k}$ is invariant under the time-reversal operator $r_T$ when restricted to the time interval $[0,T)$, we have by [10], Proposition VI.3.26, that the laws of $\{Y^k, k \ge 1\}$ are tight in $\mathbf{D}([0,T), \mathbb{R}^n)$ and any of subsequential limits of $\{Y^k, k \ge 1\}$ is the law of a continuous process. Now by [7], Proposition 3.10.4 and [7], Problem 3.25(d) on page 153, we conclude that the laws of $\{X^k, k \ge 1\}$ are tight in the space $C([0,T), \mathbb{R}^n)$. $\square$

Let $(X, \mathbf{P})$ be any of subsequential limits of $(X^n, \mathbf{P}^n_{\mu_n})$, say, along a subsequence $\{X^{n_j}, \mathbf{P}^{n_j}_{\mu_{n_j}}, j \ge 1\}$. Let $\tau_D := \inf\{t > 0 : X_t \notin D\}$ and $m_D := m|_D$. Clearly, $X_0$ has distribution $m_D$.

LEMMA 2.2. *In the above setting, for every $f \in C_c^\infty(D)$, the process $M_t^f := f(X_t) - f(X_0) - \frac{1}{2n}\int_0^t \Delta f(X_s)\, ds$ is a $\mathbf{P}$-square integrable martingale. This, in particular, implies that $\{X_t, t < \tau_D, \mathbf{P}\}$ is a Brownian motion killed upon leaving $D$, with initial distribution $m_D$ and infinitesimal generator $\frac{1}{2n}\Delta$.*

PROOF. The lemma follows easily from the invariance principle for random walks, but we sketch an argument for the sake of completeness. Recall the definition of $(X^k, \mathbf{P}_x^k, x \in D)$ and $\{\mathcal{G}_t^k, t \geq 0\}$, the $\sigma$-field generated by $X^k$. For $f \in C_c^\infty(D)$, there is a $k_0$ such that

$$\operatorname{supp}[f] \subset \{x \in D : \operatorname{dist}(x, D^c) > 2^{-2k_0}\sqrt{n}\}.$$

Thus, by (2.4), $2^{2k}\mathcal{L}_k f$ converges uniformly to $\frac{1}{2n}\Delta f$ on $D$.

Without loss of generality, we take the sample space of $X^k$ and $X$ to be the canonical space $C([0,1], \mathbb{R}^n)$. By the same argument as that in the last paragraph on page 271 of [16], we can deduce that $\{M_t^f, t \geq 0\}$ is a $\mathbf{P}$-martingale. □

LEMMA 2.3. *Let $D$ be a bounded domain in $\mathbb{R}^n$ and fix $k \geq 1$. Then for every $j \geq 1$ and $f \in L^2(D, m_k)$,*

$$(f - Q_k^{2j}f, f)_{L^2(D,m_k)} \leq j(f - Q_k^2 f, f)_{L^2(D,m_k)} \leq 2j(f - Q_k f, f)_{L^2(D,m_k)}.$$

PROOF. Note that $Q_k$ is a symmetric operator in $L^2(D, m_k)$. So for $f \in L^2(D, m_k)$,

$$(Q_k^2 f - Q_k^4 f, f)_{L^2(D,m_k)} = (f - Q_k^2 f, Q_k^2 f)_{L^2(D,m_k)}$$
$$\leq (f - Q_k^2 f, f)_{L^2(D,m_k)}.$$

Moreover,

$$(Q_k^2 g, g)_{L^2(D,m_k)} = (Q_k g, Q_k g)_{L^2(D,m_k)} \geq 0.$$

We apply the last remark to $g := f - Q_k^2 f$ to obtain

$$(Q_k^4 f - Q_k^6 f, f)_{L^2(D,m_k)} = (Q_k^2 f - Q_k^4 f, Q_k^2 f)_{L^2(D,m_k)}$$
$$= (Q_k^2 f - Q_k^4 f, f)_{L^2(D,m_k)}$$
$$\quad - (Q_k^2(f - Q_k^2 f), f - Q_k^2 f)_{L^2(D,m_k)}$$
$$\leq (Q_k^2 f - Q_k^4 f, f)_{L^2(D,m_k)}.$$

Suppose the following holds for some $j \geq 2$:

$$(Q_k^{2i}f - Q_k^{2(i+1)}f, f)_{L^2(D,m_k)} \leq (Q_k^{2(i-1)}f - Q_k^{2i}f, f)_{L^2(D,m_k)}$$

for every $i \leq j$.

Then
$$(Q_k^{2(j+1)}f - Q_k^{2(j+2)}f, f)_{L^2(D,m_k)}$$
$$= (Q_k^{2(j-1)}(Q_k^2 f) - Q_k^{2j}(Q_k^2 f), Q_k^2 f)_{L^2(D,m_k)}$$
$$\leq (Q_k^{2(j-2)}(Q_k^2 f) - Q_k^{2(j-1)}(Q_k^2 f), Q_k^2 f)_{L^2(D,m_k)}$$
$$= (Q_k^{2j}f - Q_k^{2(j+1)}f, f)_{L^2(D,m_k)}.$$

This proves by induction that, for every $i \geq 1$,
$$(Q_k^{2i}f - Q_k^{2(i+1)}f, f)_{L^2(D,m_k)} \leq (Q_k^{2(i-1)}f - Q_k^{2i}f, f)_{L^2(D,m_k)}.$$

It follows that
$$(Q_k^{2i}f - Q_k^{2(i+1)}f, f)_{L^2(D,m_k)} \leq (f - Q_k^2 f, f)_{L^2(D,m_k)} \qquad \text{for every } i \geq 1,$$

and so

(2.5)
$$(f - Q_k^{2j}f, f)_{L^2(D,m_k)} = \sum_{i=1}^{j}(Q_k^{2(i-1)}f - Q_k^{2i}f, f)_{L^2(D,m_k)}$$
$$\leq j(f - Q_k^2 f, f)_{L^2(D,m_k)}.$$

Since $Q_k$ is a symmetric operator in $L^2(D, m_k)$, we have
$$(f - Q_k^2 f, f)_{L^2(D,m_k)} = (f - Q_k f, f)_{L^2(D,m_k)} + (f - Q_k f, Q_k f)_{L^2(D,m_k)}$$
$$\leq 2(f - Q_k f, f)_{L^2(D,m_k)}.$$

This and (2.5) prove the lemma. $\square$

We will say that "$Z_t$ is a Brownian motion running at speed $1/n$" if $Z_{nt}$ is the standard Brownian motion and we will apply the same phrase to the other related process.

THEOREM 2.4. *Let $D$ be a bounded domain in $\mathbb{R}^n$ with $m(\partial D) = 0$ and satisfying the condition (1.1). Then for every $T > 0$, the laws of $\{X^k, \mathbf{P}_{m_k}^k\}$ converge weakly in $C([0,T], \mathbb{R}^n)$ to a stationary reflected Brownian motion on $D$ running at speed $1/n$ whose initial distribution is the Lebesgue measure in $D$.*

PROOF. Fix $T > 0$. Let $(X, \mathbf{P})$ be any of the subsequential limits of $(X^k, \mathbf{P}_{m_k}^k)$ in $C([0,T], \mathbb{R}^n)$, say, along $(X^{k_j}, \mathbf{P}_{m_{k_j}}^{k_j})$. Clearly, $X$ is a time-homogeneous Markov process with transition semigroup $\{P_t, t \geq 0\}$ that is symmetric in $L^2(D, dx)$. Let $\{P_t^k, t \in 2^{-k}\mathbb{Z}_+\}$ be defined by $P_t^k f(x) :=$ $\mathbf{E}_x^k[f(X_t^k)]$. For dyadic $t > 0$, say, $t = j_0/2^{2k_0}$ and $f \in C^1(\overline{D})$, we have, by Lemma 2.3 and the mean-value theorem,

$$\frac{1}{t}(f - P_t f, f)_{L^2(D, dx)}$$

$$= \frac{1}{t} \lim_{j \to \infty} (f - P_t^{k_j} f, f)_{L^2(D, m_{k_j})}$$

$$= \frac{2^{2k_0}}{j_0} \lim_{j \to \infty} (f - Q_{k_j}^{j_0 2^{2k_j - 2k_0}} f, f)_{L^2(D, m_{k_j})}$$

$$\leq \limsup_{j \to \infty} \frac{2^{2k_0}}{j_0} j_0 2^{2k_j - 2k_0} (f - Q_{k_j} f, f)_{L^2(D, m_{k_j})}$$

$$= \limsup_{j \to \infty} 2^{(2-n)k_j} \frac{1}{2n} \sum_{x \in D_{k_j}} \sum_{y \in D_{k_j} : y \leftrightarrow x} (f(x)^2 - f(x)f(y))$$

$$= \limsup_{j \to \infty} 2^{(2-n)k_j} \frac{1}{4n} \sum_{x \in D_{k_j}} \sum_{y \in D_{k_j} : y \leftrightarrow x} (f(x)^2 - f(x)f(y))$$

$$+ \limsup_{j \to \infty} 2^{(2-n)k_j} \frac{1}{4n} \sum_{y \in D_{k_j}} \sum_{x \in D_{k_j} : y \leftrightarrow x} (f(y)^2 - f(x)f(y))$$

$$= \limsup_{j \to \infty} 2^{(2-n)k_j} \frac{1}{4n} \sum_{x, y \in D_{k_j} : y \leftrightarrow x} (f(x) - f(y))^2$$

$$= \frac{1}{2n} \int_D |\nabla f(x)|^2 \, dx.$$

Let $(\mathcal{E}, \mathcal{F})$ be the Dirichlet form of $X$, or equivalently, of semigroup $\{P_t, t \geq 0\}$. That is,

$$\mathcal{F} = \left\{ f \in L^2(D, dx) : \sup_{t > 0} \frac{1}{t} (f - P_t f, f)_{L^2(D, dx)} \right.$$

$$\left. = \lim_{t \to 0} \frac{1}{t} (f - P_t f, f)_{L^2(D, dx)} < \infty \right\},$$

$$\mathcal{E}(f, f) = \sup_{t > 0} \frac{1}{t} (f - P_t f, f)_{L^2(D, dx)}$$

$$= \lim_{t \to 0} \frac{1}{t} (f - P_t f, f)_{L^2(D, dx)} \qquad \text{for } f \in \mathcal{F}.$$

Then for $f \in C^1(\overline{D})$,

$$\mathcal{E}(f, f) = \lim_{t \to 0} \frac{1}{t} (f - P_t f, f)_{L^2(D, dx)} \leq \frac{1}{2n} \int_D |\nabla f(x)|^2 \, dx.$$

This shows that $f \in \mathcal{D}(\mathcal{E})$. As $C^1(\overline{D})$ is dense in $(W^{1,2}(D), \|\cdot\|_{1,2})$, we have $W^{1,2}(D) \subset \mathcal{F}$ and

$$\mathcal{E}(f,f) \leq \frac{1}{2n} \int_D |\nabla f(x)|^2 \, dx \qquad \text{for every } f \in W^{1,2}(D).$$

This, Lemma 2.2 and Theorem 1.1 imply that $\mathcal{F} = W^{1,2}(D)$ and

$$\mathcal{E}(f,f) = \frac{1}{2n} \int_D |\nabla f(x)|^2 \, dx \qquad \text{for } f \in W^{1,2}(D).$$

We deduce then that $X$ is a stationary reflected Brownian motion on $D$ running at speed $1/n$. This proves that $X^k$ converge weakly on $C([0,T],\mathbb{R}^n)$ to the stationary reflected Brownian motion on $D$. □

REMARK 2.5. By [7], Proposition 3.10.4, under the assumptions of Theorem 2.4, the stationary laws $\mathbf{P}^k_{m_k}$ of the step processes $Y^k$ defined at the beginning of this section converge weakly in $\mathbf{D}([0,T],\mathbb{R}^n)$ to the stationary reflected Brownian motion on $D$ running at speed $1/n$, for every $T > 0$.

**3. Continuous-time reflected random walk.** Let $D$ be a bounded domain in $\mathbb{R}^n$ with $m(\partial D) = 0$ and let $D_k$ be defined as in the previous section. But in this section, $X^k$ will be the continuous time simple random walk on $D_k$, making jumps at the rate $2^{-2k}$. By definition, $X^k$ jumps to one of its nearest neighbors with equal probabilities. This process is symmetric with respect to measure $m_k$, where $m_k(x) = \frac{v_k(x)}{2n} 2^{-kn}$ for $x \in D_k$. Note that $m_k$ converge weakly to the Lebesgue measure $m$ on $D$, and recall that, for $x, y \in D_k$, we write $x \leftrightarrow y$ if $x$ and $y$ are at the distance $2^{-k}$. The Dirichlet form of $X^k$ is given by

$$\mathcal{E}^k(f,f) = \frac{1}{4n} \sum_{x,y \in D_k : x \leftrightarrow y} 2^{-(n-2)k}(f(x) - f(y))^2.$$

LEMMA 3.1. Let $D$ be a bounded domain in $\mathbb{R}^n$ with $m(\partial D) = 0$. Then for $f \in C^1(\overline{D})$,

(3.1) $$\lim_{k \to \infty} \mathcal{E}^k(f,f) = \frac{1}{2n} \int_D |\nabla f(x)|^2 m(dx).$$

PROOF. For $\delta > 0$, define $D^\delta := \{x \in \mathbb{R}^n : \text{dist}(x, \partial D) < \delta\}$. As $m(\partial D) = 0$, for every $\varepsilon > 0$, there is $\delta > 0$ such that $m(D^\delta) < \varepsilon$. Take integer $k_0 \geq 1$ large enough so that $2^{-2k_0} < \delta$. Recall that for $x \in D_k$, $v_k(x)$ denotes the degree of vertex $x$ in the graph $D_k$. Define $S_k := \{x \in D_k : v(k) = 2n\}$. Then for $k \geq k_0$, $(D \setminus D^\delta) \cap 2^{-k}\mathbb{Z}^n \subset S_k$. As

$$\mathcal{E}^k(f,f) = \frac{1}{4n} \sum_{x \in D_k} \left( \sum_{y \in D_k : y \leftrightarrow x} 2^{-(n-2)k}(f(x) - f(y))^2 \right),$$

we have by the mean-value theorem that

$$\limsup_{k\to\infty} \left| \mathcal{E}^k(f,f) - \frac{1}{2n}\int_{D\setminus D^\delta} |\nabla f(x)|^2 m(dx) \right| \leq O(m(D^\delta)).$$

Taking $\delta \downarrow 0$ yields the claim that $\lim_{k\to\infty} \mathcal{E}^k(f,f) = \frac{1}{2n}\int_D |\nabla f(x)|^2 m(dx)$ for $f \in C^1(\overline{D})$. $\square$

Let $\mathbf{P}^k_{m_k}$ denote the distribution of $\{X^k_t, t \geq 0\}$ with the initial distribution $m_k$.

LEMMA 3.2. *Assume that $D$ is a bounded domain in $\mathbb{R}^n$ with $m(\partial D) = 0$. For every $T > 0$, the laws of stationary random walks $\{X^k, \mathbf{P}^k_{m_k}, k \geq 1\}$ are tight in the space $\mathbf{D}([0,T], \overline{D})$ equipped with the Skorokhod topology.*

PROOF. For constant $T > 0$, let $r_T$ be the time-reversal operator from time $T$ for $X^k$ [see (2.1) for its definition]. Note that $\mathbf{P}^k_{m_k}$-a.s., $X^k$ is continuous at time $T$ and $\mathbf{P}^k_{m_k}$ is invariant under $r_T$. For $f \in C^1(\overline{D})$, we have by [6], Lemma 3.5 and (3.6) the following forward–backward martingale decomposition for $f(X^k_t)$. For every $T > 0$ there exists a martingale $M^{k,f}$ such that

$$(3.2) \qquad f(X^k_t) - f(X^k_0) = \tfrac{1}{2} M^{k,f}_t - \tfrac{1}{2}(M^{k,f}_{T+1} - M^{k,f}_{(T+1-t)-}) \circ r_{T+1}$$

$$\text{for } t \in [0,T].$$

Note that the symmetric jump process $X^k$ has Lévy system $(N^k(x,dy), t)$, where, for $x \in D_k$,

$$N^k(x,dy) = \frac{1}{2n} \sum_{y \in D_k : y \leftrightarrow x} 2^{-(n-2)k} 2^{nk} \delta_{\{y\}}(dy) = \frac{1}{2n} \sum_{y \in D_k : y \leftrightarrow x} 2^{2k} \delta_{\{y\}}(dy).$$

Thus,

$$\langle M^{k,f} \rangle_t = \int_0^t (f(X^k_s) - f(y))^2 N^k(X^k_s, dy)\, ds$$

$$= \frac{1}{2n} \int_0^t \sum_{y \in D_k : y \leftrightarrow X^k_s} 2^{2k} (f(X^k_s) - f(y))^2\, ds.$$

Recall that $v(k)$ is the degree of vertex $x$ in the graph $D_k$. Taking $f(x) = x_i$, $i = 1, \ldots, n$, we have, for every $k \geq 1$ and $t > s \geq 0$, with $M^{k,i} = M^{k,x_i}$,

$$\sum_{i=1}^n (\langle M^{k,i} \rangle_t - \langle M^{k,i} \rangle_s) \leq \frac{1}{2n} \int_s^t v_k(X^k_r) 2^{2k} 2^{-2k}\, dr \leq t - s.$$

This implies that the sequence $\{\sum_{i=1}^n \langle M^{k,i} \rangle_t, k \geq 1\}$ is $C$-tight in $\mathbf{D}(\mathbb{R})$ in the sense of [10], Proposition VI.3.26. Hence, by [10], Theorem VI.4.13, the laws of $\{(M^{k,1},\ldots,M^{k,n}), k \geq 1\}$ are tight in $\mathbf{D}([0,T],\mathbb{R}^n)$. As $m_k$ converges weakly to $m$, and $\mathbf{P}^k_{m_k}$ is invariant under $r_{T+1}$ for every $k \geq 1$, we conclude by [10], Theorem VI.3.21, that the laws of $\{X^k, k \geq 1\}$ are tight in $\mathbf{D}([0,T],\mathbb{R}^n)$ and, hence, on $\mathbf{D}([0,\infty),\overline{D})$. $\square$

THEOREM 3.3. *Let $D$ be a bounded domain in $\mathbb{R}^n$ with $m(\partial D) = 0$ and satisfying the condition (1.1). Then for every $T > 0$, the stationary random walks $X^k$ on $D_k$ converge weakly in the space $\mathbf{D}([0,T],\overline{D})$, as $k \to \infty$, to the stationary reflected Brownian motion on $\overline{D}$ running at speed $1/n$, whose initial distribution is the Lebesgue measure in $D$.*

PROOF. Let $(Z, \mathbf{P})$ be any of the subsequential limits of $(X^k, \mathbf{P}^k_{m_k})$, say, along $X^{k_j}$. Clearly, $Z$ is a time-homogeneous Markov process under $\mathbf{P}$ and its transition semigroup $\{P_t, t \geq 0\}$ is symmetric with respect to the Lebesgue measure $m$ on $D$. By a similar argument as that in the proof of Lemma 2.2, the process $Z$ killed upon leaving $D$ is a killed Brownian motion in $D$ with speed $1/n$. Let $(\mathcal{E}, \mathcal{F})$ be the Dirichlet form associated with $Z$, and let $\{P^k_t, t \geq 0\}$ be the transition semigroup for $X^k$. As $X^{k_j}$ converge weakly to $Z$, we have, for every $f \in C^2(\overline{D})$ and $t > 0$,

$$\lim_{j \to \infty} \frac{1}{t}(f - P^{k_j}_t f, f) = \lim_{j \to \infty} \frac{1}{t} \mathbf{E}_{\mathbf{P}^{k_j}_{m_{k_j}}}[f(X^{k_j}_0)(f(X^{k_j}_0) - f(X^{k_j}_t))]$$

$$= \frac{1}{t}\mathbf{E}_{\mathbf{P}}[f(Z_0)(f(Z_0) - f(Z_t))] = \frac{1}{t}(f - P_t f, f).$$

Thus, for $f \in C^2(\overline{D})$, by Lemma 3.1,

$$\mathcal{E}(f,f) = \sup_{t>0} \frac{1}{t}(f - P_t f, f)$$

$$= \sup_{t>0} \lim_{j \to \infty} \frac{1}{t}(f - P^{k_j}_t f, f)$$

$$\leq \liminf_{j \to \infty} \sup_{t>0} \frac{1}{t}(f - P^{k_j}_t f, f)$$

$$= \liminf_{j \to \infty} \mathcal{E}^{k_j}(f,f)$$

$$= \frac{1}{2n} \int_D |\nabla f(x)|^2 m(dx).$$

By assumption (1.1), $C^1(\overline{D})$ is dense in the Sobolev space $W^{1,2}(D)$ with respect to norm $\|\cdot\|_{1,2}$. It follows that $\mathcal{F} \supset W^{1,2}(D)$ and

$$\mathcal{E}(f,f) \leq \frac{1}{2n} \int_D |\nabla f(x)|^2 m(dx) \qquad \text{for every } f \in W^{1,2}(D).$$

Define

$$\mathcal{E}^0(f,g) = \frac{1}{2n} \int_D \nabla f(x) \cdot \nabla g(x) m(dx) \qquad \text{for } f, g \in W^{1,2}(D).$$

Note that $(\mathcal{E}^0, W^{1,2}(D))$ is the Dirichlet form for the reflected Brownian motion on $D$ running at speed $1/n$. On the other hand, as we have observed at the beginning of this proof, the process $Z$ killed upon leaving $D$ is a killed Brownian motion in $D$ with speed $1/n$. Therefore, according to Theorem 1.1, $(\mathcal{E}, \mathcal{F}) = (\mathcal{E}^0, W^{1,2}(D))$. In other words, we have shown that every subsequential limit of $X^k$ is reflected Brownian motion on $D$ with initial distribution being the Lebesgue measure on $D$ and with speed $1/n$. This shows that $X^k$ converges weakly on the space $\mathbf{D}([0,\infty), \overline{D})$ to the stationary reflected Brownian motion $X$ on $D$ running at speed $1/n$. □

**4. Examples.** All the results in the previous two sections apply to any bounded domain $D$ that satisfies the condition (1.1) and whose boundary has zero Lebesgue measure. As we noted in Section 1, bounded uniform domains have such properties. Bounded Lipschitz domains and bounded nontangentially accessible domains are uniform domains. Although the Hausdorff dimension of the Euclidean boundary of any uniform domain $D \subset \mathbb{R}^n$ is strictly less than $n$ and, thus, $\partial D$ has zero Lebesgue measure in $\mathbb{R}^n$, $\partial D$ can be highly nonrectifiable. The classical "von Koch snowflake" planar domain defined below is such an example.

To define the von Koch snowflake, start with an equilateral triangle $T_1$. Let $I$ be any of its sides. Add an equilateral triangle such that one of its sides is the middle one third of $I$ and its interior does not intersect $T_1$. There are three such triangles; let $T_2$ be the closure of the union of these three triangles and $T_1$. We proceed inductively. Suppose $I$ is one of the line segments in $\partial T_j$. Add an equilateral triangle such that one of its sides is the middle one third of $I$ and its interior does not intersect $T_j$. Let $T_{j+1}$ be the closure of the union of all such triangles and $T_j$. The snowflake $D_{\text{vK}}$ is the interior of the closure of the union of all triangles constructed in all inductive steps.

It is elementary to check that $D_{\text{vK}}$ is a nontangentially accessible domain and so a uniform domain or an $(\varepsilon, \infty)$-domain (see Section 1 for definitions). It is also well known that the Hausdorff dimension of $\partial D_{\text{vK}}$ is $\frac{\log 4}{\log 3}$. Hence, the 2-dimensional Lebesgue measure of $\partial D_{\text{vK}}$ is 0. We conclude that all results stated in the previous two sections, in particular, Theorems 2.4 and 3.3, apply to the von Koch snowflake.

Without some domain regularity conditions, the results in the previous two sections do not have to be true. Here is a counter-example. Let

$$U_k^\varepsilon = \{(x,y) \in (0,1)^2 : |x - j2^{-k}| < \varepsilon \text{ or } |y - j2^{-k}| < \varepsilon \text{ for some } j \in \mathbb{Z}\}.$$

We choose $\varepsilon_k > 0$ so that $|U_k^{\varepsilon_k}| < 2^{-k-1}$ and let $U = \bigcup_{k \geq 1} U_k^{\varepsilon_k}$. Note that $U$ is a bounded open connected set with Lebesgue area less than $1/2$. Let $D_k$ be defined as in the previous section, relative to $D = (0,1)^2$. Note that for every $k \geq 1$, $D_k \subset U$, so the sets $D_k$ defined relative to $U$ are the same as those defined relative to $D = (0,1)^2$. It follows from Theorems 2.4 and 3.3 that processes $X^k$, with discrete and continuous time, defined relative to $D_k$, converge weakly to the reflected Brownian motion in $[0,1]^2$. Since $m(U) \neq m((0,1)^2)$, the reflected Brownian motion in $U$ has a different distribution than the reflected Brownian motion in $(0,1)^2$ and it follows that the conclusions of Theorems 2.4 and 3.3 do not hold for $U$.

We would like to emphasize the fact that the approximation scheme for reflected Brownian motion developed in the next section works for *any* bounded domain.

**5. Myopic conditioning.** Throughout this section, $D \subset \mathbb{R}^n$ is a bounded connected open set and $X = \{X_t, t \geq 0, \mathbf{P}_x, x \in \mathbb{R}^n\}$ is a Brownian motion in $\mathbb{R}^n$. Let $X^D = \{X_t^D, t \geq 0, \mathbf{P}_x, x \in D\}$ be the Brownian motion in $\mathbb{R}^n$ killed upon leaving the domain $D$. For each $k \geq 1$, we define a myopic process $X^k$ as follows. For $t \in [0, 2^{-k}]$, let $X_t^k$ be $X^D$ conditioned not to leave domain $D$ by time $2^{-k}$. Suppose that $X^k$ is now defined on the time interval $[0, j2^{-k}]$. We define $X_{j2^{-k}+s}^k$ for $s \in (0, 2^{-k}]$ to be a copy of $X^D$ conditioned not to leave domain $D$ by time $2^{-k}$, starting from $X_{j2^{-k}}^k$, but otherwise independent of $\{X_t^k, t \in [0, j2^{-k}]\}$. The law of $X^k$ with $X_0^k = x$ will be denoted as $\mathbf{P}_x^k$ and the mathematical expectation under $\mathbf{P}_x^k$ will be denoted by $\mathbf{E}_x^k$.

THEOREM 5.1. *For every bounded domain $D$ in $\mathbb{R}^n$ and for every $x_0 \in D$, the processes $X^k$ under $\mathbf{P}_{x_0}^k$ converge weakly to the reflected Brownian motion $Y$ on $\overline{D}$ starting from $x_0$ in the space $C([0,1], \overline{D})$ as $k \to \infty$.*

To prove the above theorem, we introduce auxiliary processes in which pieces of conditioned Brownian paths are replaced with line segments. More precisely, we let $\{Y^k, t \in [0,1]\}$ be constructed from $\{X_{j2^{-k}}^k, j = 0, 1, \ldots, 2^k\}$ by linear interpolation over the intervals $((j-1)2^{-k}, j2^{-k})$ for $j = 1, \ldots, 2^k$. Note that $\{Y_{j2^{-k}}^k, j = 0, 1, \ldots, 2^k\}$ is a Markov chain with one-step transition probability $Q_k$, where, for $x \in D$ and $k \geq 1$, $Q_k(x, dy)$ is the distribution at time $2^{-k}$ of $X^D$ that starts from $x$ and is conditioned not to leave $D$ by time $2^{-k}$. In other words, if we let $\{P_t^D, t \geq 0\}$ denote the transition semigroup for the killed Brownian motion $X^D$, then for any nonnegative Borel function $f$ on $D$,

$$Q_k f(x) := \int_D f(y) Q_k(x, dy) = \frac{P_{2^{-k}}^D f(x)}{P_{2^{-k}}^D 1(x)}.$$

With a slight abuse of notation, the law of $Y^k$ with $Y_0^k = x$ will also be denoted as $\mathbf{P}_x^k$ and the mathematical expectation under $\mathbf{P}_x^k$ will be denoted by $\mathbf{E}_x^k$. Let $m_k(dx) := 1_D(x) P_{2^{-k}}^D 1(x)\, dx$. Observe that

(5.1) $\quad (Q_k f, g)_{L^2(D, m_k)} = (f, Q_k g)_{L^2(D, m_k)} \quad$ for $f, g \geq 0$ on $D$,

and so $m_k$ is a reversible measure for Markov chain $\{Y_{j2^{-k}}^k, j = 0, 1, \ldots, 2^k\}$. Let $m_D$ denote the Lebesgue measure on $D$ that is extended to $\mathbb{R}^n$ by letting $m_D(\mathbb{R}^n \setminus D) = 0$. It is clear that $m_k$ converge weakly to $m_D$ on $D$.

We will show that $Y^k$ converge weakly to reflected Brownian motion on $\overline{D}$ and then use this fact to establish Theorem 5.1.

LEMMA 5.2. *Suppose that either* (i) $\mu_k = m_k$ *for every* $k \geq 1$; *or* (ii) $\{\mu_k, k \geq 1\}$ *is a sequence of measures on* $D$ *with* $\sup_{k \geq 1} \mu_k(D) < \infty$ *and* $\mu_k(D \setminus K) = 0$ *for some compact subset* $K$ *of* $D$ *and all* $k \geq 1$. *Then the laws of* $\{Y^k, \mathbf{P}_{\mu_k}^k, k \geq 1\}$ *are tight in the space* $C([0,1], \mathbb{R}^n)$.

PROOF. (i) We first prove the lemma under condition (i). For nonnegative $f \in C^2(\overline{D})$, by Itô's formula,

$$Q_k f(x) = \frac{1}{P_{2^{-k}}^D 1(x)} \mathbf{E}_x [f(X_{2^{-k}}^D)]$$

$$= \frac{1}{P_{2^{-k}}^D 1(x)} \left( f(x) + \frac{1}{2} \mathbf{E}_x \left[ \int_0^{2^{-k}} \Delta f(X_s^D)\, ds \right] \right)$$

$$\geq f(x) - \frac{\|\Delta f\|_\infty}{2} 2^{-k}.$$

Fix $k \geq 1$. Let $\mathcal{G}_{j2^{-k}}^k = \sigma(Y_{i2^{-k}}^k, i \leq j)$. For nonnegative $f \in C^2(\overline{D})$, let $A_f := \frac{\|\Delta f\|_\infty}{2}$. Then we see from the above that

$$\mathbf{E}_{m_k}^k [f(Y_{(j+1)2^{-k}}^k) + A_f(j+1)2^{-k} | \mathcal{G}_{j2^{-k}}^k] = Q_k f(Y_{j2^{-k}}^k) + A_f(j+1)2^{-k}$$

$$\geq f(Y_{j2^{-k}}^k) + A_f j 2^{-k}.$$

In other words, $\{f(Y_{j2^{-k}}^k) + A_f j 2^{-k}, \mathcal{G}_{j2^{-k}}^k\}_{j=0,1,\ldots,2^k}$ is a nonnegative $\mathbf{P}_{m_k}^k$-submartingale. Moreover, for every $\varepsilon > 0$,

$$\lim_{k \to \infty} \sum_{j=1}^{2^k} \mathbf{P}_{m_k}^k (|Y_{j2^{-k}}^k - Y_{(j-1)2^{-k}}^k| > \varepsilon)$$

$$\leq \lim_{k \to \infty} 2^k \int_D \mathbf{P}_x(|X_{2^{-k}} - X_0| > \varepsilon \text{ and } 2^{-k} < \tau_D)\, dx$$

$$\leq \lim_{k \to \infty} 2^k \int_D \mathbf{P}_x(|X_{2^{-k}} - X_0| > \varepsilon)\, dx$$

$$= 0.$$

Thus, by [16], Theorem 1.4.11, the laws of $\{Y^k, \mathbf{P}^k_{m_k}, k \geq 1\}$ are tight in $C([0,1], \mathbb{R}^n)$.

(ii) Now assume that $\{\mu_k, k \geq 1\}$ is a sequence of measures on $D$ satisfying condition (ii). Since every compact set has a finite covering by open balls, we may and do assume, without loss of generality, that $K \subset \overline{B(x_0, r_0)} \subset D$. Define $\delta_0 := \operatorname{dist}(B(x_0, r_0), D^c)/2$ and $B := B(x_0, r_0 + \delta_0)$. Recall that $X$ is Brownian motion in $\mathbb{R}^n$. By Lemma II.1.2 in Stroock [14], there is a constant $c_0 > 0$ so that

$$(5.2) \qquad \mathbf{P}\left(\sup_{s \leq t} |X_s - X_0| > r\right) \leq c_0 \exp\left(-\frac{r}{c_0 t}\right).$$

This, in particular, implies that

$$(5.3) \qquad P_t^D \mathbf{1}(x) \geq 1 - c_0 \exp\left(-\frac{\delta_0}{c_0 t}\right) \qquad \text{for every } x \in B.$$

Let $P_t^B$ denote the semigroup for Brownian motion killed upon exiting $B$.

Without loss of generality, we take the sample space of $Y^k$ and $X$ to be the canonical space $C([0,1], \mathbb{R}^n)$. For $\omega \in C([0,1], \mathbb{R}^n)$ and $\rho > 0$, we define the oscillation of $\omega$ over time interval $[s_0, t_0]$ by

$$\operatorname{osc}_\rho[s_0, t_0](\omega) := \sup_{s,t \in [s_0, t_0]:\, |t-s| \leq \rho} |\omega(s) - \omega(t)|.$$

We also define

$$\tau_0 := \inf\{t \geq 0 : \omega(t) \notin B\}.$$

Sometimes we will add a superscript to make the underlying process explicit, for example, we may write $\mathbf{P}(\operatorname{osc}^X_\rho[s_0, t_0] > \varepsilon)$.

Since all $\mu_k$'s are supported in a compact set and have uniformly bounded mass, standard theorems (see, e.g., [16], Theorem 1.3.1) show that the laws of $\{Y^k, \mathbf{P}^k_{\mu_k}, k \geq 1\}$ are tight in the space $C([0,1], \mathbb{R}^n)$ if and only if, for every $\varepsilon > 0$,

$$\limsup_{\rho \downarrow 0} \sup_{k \geq 1} \mathbf{P}^k_{\mu_k}(\operatorname{osc}_\rho[0,1] > \varepsilon) = 0.$$

So it suffices to show that, for every $\varepsilon > 0$ and $\delta > 0$, there is $\rho > 0$ and $N \geq 1$ such that

$$(5.4) \qquad \mathbf{P}^k_{\mu_k}(\operatorname{osc}_\rho[0,1] > \varepsilon) < \delta \qquad \text{for every } k \geq N.$$

Fix an arbitrarily small $\delta > 0$, and let $t_0 = j_0 2^{k_0}$ be a dyadic rational in $(0,1)$ such that $c_0 \exp(-\frac{\delta_0}{c_0 t_0}) \sup_{k \geq 1} \mu_k(D) < \delta/4$. Define

$$c_k := \left(1 - c_0 \exp\left(-\frac{2^k \delta_0}{c_0}\right)\right)^{-2^k t_0},$$

and note that $c_k \to 1$ as $k \to \infty$. We have

$$\mathbf{P}^k_{\mu_k}(\mathrm{osc}_\rho[0,1] > \varepsilon) \leq \mathbf{P}^k_{\mu_k}(\tau_0 \leq t_0) + \mathbf{P}^k_{\mu_k}(\mathrm{osc}_\rho[0,1] > \varepsilon \text{ and } \tau_0 > t_0).$$

Observe that $B$ is convex. For $k \geq k_0$, it follows from the definition of the process $Y^k$ and (5.3) that

$$\mathbf{P}^k_{\mu_k}(\tau_0 \leq t_0)$$
$$= \sum_{j=1}^{2^k t_0} \mathbf{P}^k_{\mu_k}(\tau_0 \in ((j-1)2^{-k}, j2^{-k}])$$
$$= \sum_{j=1}^{2^k t_0} \mathbf{E}_{\mu_k}\left[\left(\prod_{i=0}^{j-1} \frac{1}{P^D_{2^{-k}}1(X^D_{j2^{-k}})}\right); \right.$$
$$\left. X^D_{i2^{-k}} \in B \text{ for } i = 0, \ldots, j-1 \text{ but } X^D_{j2^{-k}} \notin B\right]$$

(5.5)
$$\leq c_k \sum_{j=1}^{2^k t_0} \mathbf{P}_{\mu_k}(X^D_{i2^{-k}} \in B \text{ for } i=0,\ldots,j-1 \text{ but } X^D_{j2^{-k}} \notin B)$$
$$\leq c_k \mathbf{P}_{\mu_k}(\tau^B_0 \leq t_0)$$
$$\leq c_k c_0 \exp\left(-\frac{\delta_0}{c_0 t_0}\right)\mu_k(D)$$
$$\leq c_k \delta/4.$$

On the other hand, by the definition of $Y^k$ and (5.3) again, conditioned on $\{\tau_0 > t_0\}$, the law of $\mathbf{P}^k_{\mu_k}$ restricted to $[0, t_0]$ is dominated by $c_k$ times that of linear interpolation of $X^D$ at times $j2^{-k}$, with initial distribution $\mu_k$. Thus, for any $\varepsilon > 0$, one can make $t_0 > 0$ smaller, if necessary, so that for all sufficiently large $k$,

(5.6)
$$\mathbf{P}^k_{\mu_k}(\tau_0 > t_0 \text{ and } \mathrm{osc}_\rho[0, t_0] > \varepsilon/2)$$
$$\leq c_k \mathbf{P}_{\mu_k}(\tau^X_0 > t_0 \text{ and } \mathrm{osc}^X_\rho[0, t_0] > \varepsilon/2)$$
$$\leq c_k \mu_k(D)\mathbf{P}_0(\mathrm{osc}^X_\rho[0, t_0] > \varepsilon/2), \leq \delta/4,$$

where $\mathbf{P}_0$ denotes the law of Brownian motion $X$ starting from the origin.

For each fixed $y \in \overline{B(x_0, r_0)}$, let $x \mapsto \psi_k(y, x)$ be the density function for the distribution of $X^D_{t_0}$ under $\mathbf{P}_y$ restricted on the event that $\{X^D_{j2^{-k}} \in B \text{ for } j = 1, 2, \ldots, 2^k t_0\}$. Clearly, $x \mapsto \psi_k(y, x)$ is a bounded function on $D$ that vanishes outside $B$, and as $k \uparrow \infty$, $\psi_k(y, x)$ decrease to $\phi(y, x)$, the probability density function for the killed Brownian motion in $B$ at time $t_0$ starting from $y \in \overline{B(x_0, r_0)}$. There is a constant $a_1 > 0$ such that

$$\sup_{y \in \overline{B(x_0, r_0)}} \psi_k(y, x) \, dx \le a_1 P^D_{t_0} 1(x) \, dx \le a_1 m_k(dx)$$

for every $k \ge 1$ such that $2^{-k} \le t_0$. As $\sup_{k \ge 1} \mu_k(D) < \infty$, there is a constant $a_2 > 0$, independent of $k \ge 1$, such that the distribution of $X^D_{t_0}$ under $\mathbf{P}_{\mu_k}$ restricted on the event that $\{X^D_{j2^{-k}} \in B \text{ for } j = 1, 2, \ldots, 2^k t_0\}$ is dominated by $a_2 m_k(dx)$. We obtain,

$$\mathbf{P}^k_{\mu_k}(\tau_0 > t_0 \text{ and } \operatorname{osc}_\rho[t_0, 1] > \varepsilon/2)$$

(5.7)
$$= \mathbf{E}_{\mu_k}\left[\left(\prod_{j=0}^{2^k t_0} \frac{1}{P^D_{2^{-k}} 1(X^D_{j2^{-k}})}\right) \mathbf{P}^k_{X^D_{t_0}}(\operatorname{osc}_\rho[0, 1 - t_0] > \varepsilon/2);\right.$$
$$\left. X^D_{j2^{-k}} \in B \text{ for } j = 1, 2, \ldots, 2^k t_0\right]$$

$$\le c_k^{(2^k t_0 + 1)/(2^k t_0)} a_2 \mathbf{P}^k_{m_k}(\operatorname{osc}_\rho[0, 1] > \varepsilon/2).$$

We have already proved in part (i) of the proof that the laws of $\{Y^k, \mathbf{P}^k_{m_k}, k \ge 1\}$ are tight in $C([0, 1], \mathbb{R}^n)$, so there is $N \ge 1$ such that

(5.8) $\quad\quad \mathbf{P}^k_{m_k}(\operatorname{osc}_\rho[0, 1] > \varepsilon/2) < \delta/(4a_2) \quad$ for every $k \ge N$.

Combining (5.5)–(5.8), we obtain for large $k$,

$$\mathbf{P}^k_{\mu_k}(\operatorname{osc}_\rho[0, 1] > \varepsilon)$$
$$\le \mathbf{P}^k_{\mu_k}(\tau_0 \le t_0) + \mathbf{P}^k_{\mu_k}(\operatorname{osc}_\rho[0, t_0] > \varepsilon/2 \text{ and } \tau_0 > t_0)$$
$$+ \mathbf{P}^k_{\mu_k}(\operatorname{osc}_\rho[t_0, 1] > \varepsilon/2 \text{ and } \tau_0 > t_0)$$
$$\le c_k \delta/4 + \delta/4 + c_k^{(2^k t_0 + 1)/(2^k t_0)} a_2 \delta/(4a_2).$$

This proves (5.4) because $\lim_{k \to \infty} c_k = 1$. □

Consider a sequence of finite measures $\mu_k$ on $D$, $k \ge 1$, that satisfies the assumptions of Lemma 5.2(i) or (ii) and converges weakly to a finite measure $\mu$. By Lemma 5.2, the laws of $\{Y^k, \mathbf{P}^k_{\mu_k}, k \ge 1\}$ are tight on $C([0, 1], \mathbb{R}^n)$. Let $(Y, \mathbf{P})$ be any of its subsequential limit, say, along a subsequence $\{Y^{n_j}, \mathbf{P}^{n_j}_{\mu_{n_j}}, j \ge 1\}$. Clearly $Y_0$ has distribution $\mu$.

LEMMA 5.3. *In the above setting, for every $f \in C_c^\infty(D)$, $M^f_t := f(Y_t) - f(Y_0) - \frac{1}{2} \int_0^t \Delta f(Y_s) \, ds$ is a $\mathbf{P}$-square integrable martingale. This in particular implies that $\{Y_t, t < \tau_D, \mathbf{P}\}$, with $\tau_D := \inf\{t > 0 : Y_t \notin D\}$, is the killed Brownian motion in $D$ with initial distribution $\mu$.*

PROOF.　Recall the definition of $(Y^k, \mathbf{P}_x^k, x \in D)$ and denote by $\{\mathcal{G}_t^k, t \geq 0\}$ the $\sigma$-field generated by $Y^k$. For $f \in C_c^\infty(D)$, define

$$\mathcal{L}_k f(x) = \int_D (f(y) - f(x)) Q_k(x, dy).$$

Then $\{f(Y_{j2^{-k}}^k) - \sum_{i=1}^j \mathcal{L}_k f(Y_{i2^{-k}}^k), \mathcal{G}_{j2^{-k}}^k, j = 0, 1, \ldots, 2^k, \mathbf{P}_x^k\}$ is a martingale for every $f \in C_c^\infty(D)$. For $f \in C_c^\infty(D)$, using the Taylor expansion, we have

$$2^k \mathcal{L}_k f(x)$$
$$= 2^k \int_D \bigg( \nabla f(x)(y-x)$$
$$\qquad + \frac{1}{2} \sum_{i,j=1}^n \frac{\partial^2 f(x)}{\partial x_i \partial x_j}(y_i - x_i)(y_j - x_j) + O(1)|y-x|^3 \bigg) Q_k(x, dy)$$
$$= \frac{2^k}{P_{2^{-k}}^D 1(x)} \bigg( \sum_{i=1}^n \frac{\partial f(x)}{\partial x_i} \mathbf{E}_x[X_{2^{-k}}^i - X_0^i; 2^{-k} < \tau_D]$$
$$\qquad + \frac{1}{2} \sum_{i,j=1}^n \frac{\partial^2 f(x)}{\partial x_i \partial x_j} \mathbf{E}_x[(X_{2^{-k}}^i - X_0^i)(X_{2^{-k}}^j - X_0^j); 2^{-k} < \tau_D]$$
$$\qquad + O(1) \mathbf{E}_x[|X_{2^{-k}} - X_0|^3; 2^{-k} < \tau_D] \bigg).$$

This converges uniformly to $\frac{1}{2}\Delta f(x)$ on $D$ as $k \to \infty$.

Without loss of generality, we take the sample space of $Y^k$ and $Y$ to be the canonical space $C([0,1], \mathbb{R}^n)$. Then by the same argument as that in the last paragraph on page 271 of [16], we can deduce that $\{M_t^f, t \geq 0\}$ is a $\mathbf{P}$-martingale.　□

LEMMA 5.4.　*Let $D$ be a bounded domain in $\mathbb{R}^n$ and fix $k \geq 1$. Then for every $j \geq 1$ and $f \in L^2(D, m_k)$,*

$$(f - Q_k^j f, f)_{L^2(D, m_k)} \leq j(f - Q_k f, f)_{L^2(D, m_k)}.$$

PROOF.　The idea of the proof is similar to that for Lemma 2.3. Note that $Q_k$ is a symmetric operator in $L^2(D, m_k)$. So for $f \in L^2(D, m_k)$,

$$(Q_k f - Q_k^2 f, f)_{L^2(D, m_k)} = (f - Q_k f, Q_k f)_{L^2(D, m_k)} \leq (f - Q_k f, f)_{L^2(D, m_k)}.$$

We have

$$(Q_k g, g)_{L^2(D, m_k)} = (P_{2^{-k}}^D g, g)_{L^2(D, dx)} = \int_D P_{2^{-k-1}}^D g(x)^2 \, dx \geq 0.$$

We apply the above observation to $g := f - Q_k f$ to obtain

$$(Q_k^2 f - Q_k^3 f, f)_{L^2(D,m_k)}$$
$$= (Q_k f - Q_k^2 f, Q_k f)_{L^2(D,m_k)}$$
$$= (Q_k f - Q_k^2 f, f)_{L^2(D,m_k)} - (Q_k(f - Q_k f), f - Q_k f)_{L^2(D,m_k)}$$
$$\leq (Q_k f - Q_k^2 f, f)_{L^2(D,m_k)}.$$

Suppose the following holds for $j \geq 2$:

$$(Q_k^i f - Q_k^{i+1} f, f)_{L^2(D,m_k)} \leq (Q_k^{i-1} f - Q_k^i f, f)_{L^2(D,m_k)} \qquad \text{for every } i \leq j.$$

Then

$$(Q_k^{j+1} f - Q_k^{j+2} f, f)_{L^2(D,m_k)} = (Q_k^{j-1}(Q_k f) - Q_k^j(Q_k f), Q_k f)_{L^2(D,m_k)}$$
$$\leq (Q_k^{j-2}(Q_k f) - Q_k^{j-1}(Q_k f), Q_k f)_{L^2(D,m_k)}$$
$$= (Q_k^j f - Q_k^{j+1} f, f)_{L^2(D,m_k)}.$$

This proves by induction that, for every $i \geq 1$,

$$(Q_k^i f - Q_k^{i+1} f, f)_{L^2(D,m_k)} \leq (Q_k^{i-1} f - Q_k^i f, f)_{L^2(D,m_k)}.$$

It follows that

$$(Q_k^i f - Q_k^{i+1} f, f)_{L^2(D,m_k)} \leq (f - Q_k f, f)_{L^2(D,m_k)} \qquad \text{for every } i \geq 1,$$

and so

$$(f - Q_k^j f, f)_{L^2(D,m_k)} = \sum_{i=1}^j (Q_k^{i-1} f - Q_k^i f, f)_{L^2(D,m_k)}$$
$$\leq j(f - Q_k f, f)_{L^2(D,m_k)},$$

which proves the lemma. $\square$

THEOREM 5.5. *Let $D$ be a bounded domain in $\mathbb{R}^n$. The processes $Y^k$ under $\mathbf{P}_{m_k}^k$ converge weakly to the stationary reflected Brownian motion $Y$ on $\overline{D}$ in the space $C([0,1], \mathbb{R}^n)$ as $k \to \infty$, where $Y_0$ is distributed according to the Lebesgue measure on $D$.*

PROOF. Let $(Y, \mathbf{P})$ be any of the subsequential limits of $(Y^k, \mathbf{P}_{m_k}^k)$ in $C([0,1], \mathbb{R}^n)$, say, along $(Y^{k_j}, \mathbf{P}_{m_{k_j}}^{k_j})$. Clearly, $Y$ is a time-homogeneous Markov process with transition semigroup $\{P_t, t \geq 0\}$ that is symmetric in $L^2(\overline{D}, m_D) = L^2(D, dx)$. Let $\{P_t^k, t \in 2^{-k}\mathbb{Z}_+\}$ be defined by $P_t^k f(x) := \mathbf{E}_x^k[f(Y_t^k)]$. As $(Y^{k_j}, \mathbf{P}^{k_j}_{m_{k_j}})$ converges weakly to $(Y, \mathbf{P})$ in $C([0,1], \mathbb{R}^n)$, we have for dyadic rational $t > 0$, say, $t = j_0/2^{k_0}$ and $\phi \in C(\overline{D})$,

$$\lim_{j \to \infty} (P^{k_j}_t \phi, \phi)_{L^2(D, m_{k_j})} = \lim_{j \to \infty} \mathbf{E}^{k_j}_{m_{k_j}}[\phi(Y^{k_j}_0)\phi(Y^{k_j}_t)]$$

(5.9)
$$= \mathbf{E}[\phi(Y_0)\phi(Y_t)] = (P_t \phi, \phi)_{L^2(D,m)}.$$

Recall from (5.1) that $m_k$ is a reversible measure for operator $P^k_t$ for $t \in 2^{-k}\mathbb{Z}_+$ when $k$ is sufficiently large. For $t \in 2^{-k}\mathbb{Z}_+$, the operator $P^k_t$ has density function $p^k(t, x, y)$ with respect to the reversible measure $m_k$. The density function $p^k(t, x, y)$ is symmetric in $(x, y)$. Note that $P^k_t 1 = 1$ and so

$$\int_D p^k(t, x, y) m_k(dy) = 1 \qquad \text{for every } x \in D.$$

For a dyadic rational $t \in (0, 1]$, $g \in L^2(D, m_k)$ and large $k$,

$$(g - P^k_t g, g)_{L^2(D, m_k)}$$

$$= \int_D (g(x) - P^k_t g(x))g(x) m_k(dx)$$

$$= \int_D \left( g(x) - \int_D p^k(t, x, y) g(y) m_k(dy) \right) g(x) m_k(dx)$$

$$= \int_{D \times D} (g(x)^2 - g(x)g(y)) p^k(t, x, y) m_k(dx) m_k(dy)$$

(5.10)
$$= \tfrac{1}{2} \int_{D \times D} (g(x)^2 - g(x)g(y)) p^k(t, x, y) m_k(dx) m_k(dy)$$

$$+ \tfrac{1}{2} \int_{D \times D} (g(y)^2 - g(x)g(y)) p^k(t, x, y) m_k(dx) m_k(dy)$$

$$= \tfrac{1}{2} \int_{D \times D} (g(x) - g(y))^2 p^k(t, x, y) m_k(dx) m_k(dy)$$

$$= \tfrac{1}{2} \mathbf{E}^k_{m_k}[(g(Y^k_0) - g(Y^k_t))^2] \geq 0.$$

For later reference, we record the following consequence of the last formula,

$$(g - Q_k g, g)_{L^2(D, m_k)}$$

(5.11)
$$= (g - P^k_{2^{-k}} g, g)_{L^2(D, m_k)} = \tfrac{1}{2} \mathbf{E}^k_{m_k}[(g(Y^k_0) - g(Y^k_{2^{-k}}))^2]$$

$$= \mathbf{E}_m[(g(X_0) - g(X_{2^{-k}}))^2; 2^{-k} < \tau_D],$$

where $X$ is Brownian motion in $\mathbb{R}^n$ and $\tau_D := \inf\{t > 0 : X_t \notin D\}$. By (5.10), we have for dyadic $t > 0$ and large $k \geq 1$,

$$0 \leq (g - P^k_{2t} g, g)_{L^2(D, m_k)} = (g, g)_{L^2(D, m_k)} - (P^k_{2t} g, g)_{L^2(D, m_k)}$$

$$= (g, g)_{L^2(D, m_k)} - (P^k_t g, P^k_t g)_{L^2(D, m_k)},$$

and thus,
$$\|P_t^k g\|_{L^2(D,m_k)}^2 \leq \|g\|_{L^2(D,m_k)}^2 \leq \|g\|_{L^2(D,m)}^2.$$

Every $f \in L^2(D,m)$ can be approximated in $L^2(D,m)$ by continuous functions on $\overline{D}$, so we have from above and (5.9) that

$$\lim_{j \to \infty} (P_t^{k_j} f, f)_{L^2(D,m_{k_j})} = (P_t f, f)_{L^2(D,m)} \qquad \text{for every } f \in L^2(D,m).$$
(5.12)

Let $Z$ be the reflected Brownian motion on $\overline{D}$, obtained from quasi-continuous projection on $\overline{D}$ of the reflected Brownian motion $Z^*$ on the Martin–Kuramochi compactification $D^*$ of $D$. Recall that $Z$ behaves like Brownian motion in $D$ before $Z$ hits the boundary $\partial D$. For $f \in W^{1,2}(D)$, it admits a quasi-continuous version on $D^*$ and, therefore, $f(Z_t) = f(Z_t^*)$ is well defined. Moreover, $f(Z_t)$ has the following Fukushima's decomposition (see [8] or [9]):

$$f(Z_t) - f(Z_0) = M_t^f + N_t^f \qquad \text{for every } t \geq 0,$$

where $M^f$ is a continuous martingale additive functional of $Z^*$ with quadratic variation process $\langle M^f \rangle_t = \int_0^t |\nabla f(Z_s)|^2 \, ds$ and $N^f$ is a continuous additive functional of $Z^*$ having zero energy. In particular, we have

$$\lim_{t \to 0} \frac{1}{t} \mathbf{E}_m[(f(Z_t) - f(Z_0))^2] = \lim_{t \to 0} \frac{1}{t} \mathbf{E}_m[(M_t^f)^2] = \int_D |\nabla f(x)|^2 \, dx.$$

These observations, (5.11)–(5.12) and Lemma 5.4 imply that for $f \in W^{1,2}(D)$ and dyadic $t > 0$,

$$\frac{1}{t}(f - P_t f, f)_{L^2(D,dx)} = \frac{1}{t} \lim_{j \to \infty} (f - P_t^{k_j} f, f)_{L^2(D,m_{k_j})}$$
$$= \frac{2^{k_0}}{j_0} \lim_{j \to \infty} (f - Q_{k_j}^{j_0 2^{k_j - k_0}} f, f)_{L^2(D,m_{k_j})}$$
$$\leq \limsup_{j \to \infty} \frac{2^{k_0}}{j_0} j_0 2^{k_j - k_0} (f - Q_{k_j} f, f)_{L^2(D,m_{k_j})}$$
$$= \limsup_{j \to \infty} 2^{k_j - 1} \mathbf{E}_m[(f(Z_{2^{-k_j}}) - f(Z_0))^2; 2^{-k_j} < \tau_D]$$
$$\leq \limsup_{j \to \infty} 2^{k_j - 1} \mathbf{E}_m[(f(Z_{2^{-k_j}}) - f(Z_0))^2]$$
$$= \frac{1}{2} \int_D |\nabla f(x)|^2 \, dx.$$

Let $(\mathcal{E}, \mathcal{F})$ be the Dirichlet form of $Y$ in $L^2(\overline{D}, m_D) = L^2(D, dx)$. Then for $f \in W^{1,2}(D)$,

$$\mathcal{E}(f,f) = \lim_{t \to 0} \frac{1}{t}(f - P_t f, f)_{L^2(D,dx)} \leq \frac{1}{2} \int_D |\nabla f(x)|^2 \, dx.$$

This shows that $f \in \mathcal{F}$. So we have $W^{1,2}(D) \subset \mathcal{F}$ and

$$\mathcal{E}(f,f) \leq \tfrac{1}{2} \int_D |\nabla f(x)|^2 \, dx \qquad \text{for every } f \in W^{1,2}(D).$$

If follows from Lemma 5.3 that $\{Y_t, t < \tau_D\}$ is a Brownian motion in $D$ with initial distribution $m_D$. Therefore, we have by Theorem 1.1 that $\mathcal{F} = W^{1,2}(D)$ and

$$\mathcal{E}(f,f) = \tfrac{1}{2} \int_D |\nabla f(x)|^2 \, dx \qquad \text{for } f \in W^{1,2}(D).$$

We deduce that $Y$ is a stationary reflected Brownian motion on $D$. This proves that $Y^k$ converges weakly on $C([0,1], \mathbb{R}^n)$ to the stationary reflected Brownian motion on $D$. □

THEOREM 5.6. *Let $D$ be a bounded domain $D$ in $\mathbb{R}^n$. Let $\{\mu_k, k \geq 1\}$ be a sequence of measures on $D$ that converges weakly to $\mu$ with $\sup_{k \geq 1} \mu_k(D) < \infty$ and that there is a compact set $K \subset D$ such that $\mu_k(D \setminus K) = 0$ for every $k \geq 1$. Then the processes $Y^k$ under $\mathbf{P}^k_{\mu_k}$ converge weakly to the reflected Brownian motion $Y$ on $\overline{D}$ with initial distribution $\mu$ in the space $C([0,1], \mathbb{R}^n)$ as $k \to \infty$.*

PROOF. Using a partition of $K$ with a finite covering of open balls if necessary, we may assume that $K \subset \overline{B(x_0, r_0)} \subset D$. Let $\delta_0 = \text{dist}(B(x_0, r_0), D^c)/2$ and define $B := B(x_0, r_0 + \delta_0)$. So the distance between ball $B$ and $D^c$ is $\delta_0$. Define

$$\tau_B := \inf\{t \geq 0 : Y_t \notin B\} \quad \text{and} \quad \tau_B^k := \inf\{t \geq 0 : Y_t^k \notin B\}.$$

According to Lemma 5.2, the laws of $\{(Y^k, \mathbf{P}^k_{x_0}), k \geq 1\}$ are tight in the space $C([0,1], \mathbb{R}^n)$. Let $(Y, \mathbf{P})$ be any of the subsequential limits of $(Y^k, \mathbf{P}^k_{\mu_k})$, say, along $(Y^{k_j}, \mathbf{P}^{k_j}_{\mu_{k_j}})$. It suffices to show that the finite-dimensional distributions of $(Y, \mathbf{P})$ are the same as those of reflected Brownian motion with initial distribution $\mu$. For this purpose, take $0 < t_1 < t_2 < \cdots < t_N$ and nonnegative $f_i \in C_b(\mathbb{R}^n)$ for $i = 1, 2, \ldots, N$. We know from Lemma 5.3 that, under $\mathbf{P}$, $\{Y_t, t < \tau_D\}$ is the killed Brownian motion in $D$ with initial distribution $\mu$. Thus, for every $\varepsilon > 0$, there is a dyadic rational $t_0 = j_0 2^{-k_0} \in (0, t_1/2)$ such that

$$\mathbf{P}(\tau_B \leq t_0) < \varepsilon.$$

Since $(Y^{k_j}, \mathbf{P}^{k_j}_{\mu_{k_j}})$ converge to $(Y, \mathbf{P})$ weakly, we have

(5.13) $$\limsup_{j \to \infty} \mathbf{P}^{k_j}_{\mu_{k_j}}(\tau_B^k \leq t_0) \leq \mathbf{P}(\tau_B \leq t_0) < \varepsilon.$$

By (5.3), there is a constant $c_0 > 0$ such that

$$P_t^D 1(x) \geq 1 - c_0 \exp\left(-\frac{\delta_0}{c_0 t}\right) \quad \text{for every } x \in B.$$

Let $P_t^D$ denote the semigroup for Brownian motion killed upon exiting $D$ and $X^B$ be the Brownian motion $X$ killed upon leaving ball $B$. The ball $B$ is convex so for any bounded nonnegative function $f$ on $\overline{D}$ and $k > k_0$,

$$\mathbf{E}_{\mu_k}^k[f(Y_{t_0}^k); \tau_B > t_0]$$
$$= \mathbf{E}_{\mu_k}^k[f(Y_{t_0}^k); Y_{j2^{-k}}^k \in B \text{ for } j = 1, 2, \ldots, 2^k t_0]$$

(5.14)
$$= \mathbf{E}_{\mu_k}\left[\left(\prod_{j=0}^{2^k t_0 - 2} \frac{1}{P_{2^{-k}}^D 1(X_{j2^{-k}}^D)}\right) \frac{P_{2^{-k}}^D f(X_{t_0 - 2^{-k}}^D)}{P_{2^{-k}}^D 1(X_{t_0 - 2^{-k}}^D)}; \right.$$
$$\left. X_{j2^{-k}}^D \in B \text{ for } j = 1, 2, \ldots, 2^k t_0 - 1\right]$$
$$\leq c_k \mathbf{E}_{\mu_k}[f(X_{t_0}^D); X_{j2^{-k}}^D \in B \text{ for } j = 1, 2, \ldots, 2^k t_0],$$

where

$$c_k := \left(1 - c_0 \exp\left(-\frac{2^k \delta_0}{2c_0}\right)\right)^{-2^k t_0} \to 1 \quad \text{as } k \to \infty.$$

For each fixed $y \in \overline{B(x_0, r_0)}$, let $x \mapsto \psi_k(y, x)$ be the density function for the distribution of $X_{t_0}^D$ under $\mathbf{P}_y$ restricted on the event that $\{X_{j2^{-k}}^D \in B \text{ for } j = 1, 2, \ldots, 2^k t_0\}$. Clearly, $x \mapsto \psi_k(y, x)$ is a bounded function on $D$ that vanishes outside $B$, and as $k \uparrow \infty$, $\psi_k(y, x)$ decrease to $p^B(t_0, y, x)$, the probability density function for the killed Brownian motion in $B$ at time $t_0$ starting from $y$. We clearly have

(5.15) $$\mathbf{E}_{\mu_k}^k[f(Y_{t_0}^k); \tau_B^k > t_0] \geq \mathbf{E}_{\mu_k}[f(X_{t_0}^B)].$$

Let

$$\phi_k(x) := \int_D p^B(t_0, x, y)\mu_k(dy) \quad \text{and} \quad \phi(x) := \int_D p^B(t_0, x, y)\mu(dy).$$

and denote $(\phi m)(dx) := \phi(x)m(dx)$ and $(\phi_j m)(dx) := \phi_j(x)m(dx)$. By (5.13)–(5.15), we have, for each $j \geq j_0$,

$$\mathbf{E}_{\phi m}^{k_j}\left[\prod_{i=1}^N f_i(Y_{t_i - t_0}^{k_j})\right] - \varepsilon \prod_{i=1}^N \|f_i\|_\infty \leq \mathbf{E}_{\mu_{k_j}}^{k_j}\left[\prod_{i=1}^N f_i(Y_{t_i}^{k_j})\right]$$
$$\leq c_{k_j} \mathbf{E}_{\phi_{k_{j_0}} m}^{k_j}\left[\prod_{i=1}^N f_i(Y_{t_i - t_0}^{k_j})\right] + \varepsilon \prod_{i=1}^N \|f_i\|_\infty.$$

Let $Z$ be the reflected Brownian motion on $D$ with initial distribution $\mu$. Since $\phi$ and $\phi_{j_0}$ are bounded continuous functions with compact support $\overline{B}$, it follows from Theorem 5.5 that

$$\lim_{j \to \infty} \mathbf{E}^{k_j}_{\phi m}\left[\prod_{i=1}^{N} f_i(Y^{k_j}_{t_i-t_0})\right] = \mathbf{E}_{\phi m}\left[\prod_{i=1}^{N} f_i(Z_{t_i-t_0})\right]$$

and

$$\lim_{j \to \infty} \mathbf{E}^{k_j}_{\phi_{j_0} m}\left[\prod_{i=1}^{N} f_i(Y^{k_j}_{t_i-t_0})\right] = \mathbf{E}_{\phi_{j_0} m}\left[\prod_{i=1}^{N} f_i(Z_{t_i-t_0})\right].$$

Taking $j \to \infty$, we have

$$\mathbf{E}_{\phi m}\left[\prod_{i=1}^{N} f_i(Z_{t_i-t_0})\right] - \varepsilon \prod_{i=1}^{N} \|f_i\|_\infty \leq \liminf_{j \to \infty} \mathbf{E}^{k_j}_{x_0}\left[\prod_{i=1}^{N} f_i(Y^{k_j}_{t_i})\right]$$

$$\leq \limsup_{j \to \infty} \mathbf{E}^{k_j}_{x_0}\left[\prod_{i=1}^{N} f_i(Y^{k_j}_{t_i})\right]$$

$$\leq \mathbf{E}_{\phi_{j_0} m}\left[\prod_{i=1}^{N} f_i(Z_{t_i-t_0})\right] + \varepsilon \prod_{i=1}^{N} \|f_i\|_\infty.$$

Since $\phi_{j_0}$ converges to $\phi$ boundedly as $j_0 \to \infty$, we have

$$\mathbf{E}_{\phi m}\left[\prod_{i=1}^{N} f_i(Z_{t_i-t_0})\right] - \varepsilon \prod_{i=1}^{N} \|f_i\|_\infty \leq \liminf_{j \to \infty} \mathbf{E}^{k_j}_{x_0}\left[\prod_{i=1}^{N} f_i(Y^{k_j}_{t_i})\right]$$

$$\leq \limsup_{j \to \infty} \mathbf{E}^{k_j}_{x_0}\left[\prod_{i=1}^{N} f_i(Y^{k_j}_{t_i})\right]$$

$$\leq \mathbf{E}_{\phi m}\left[\prod_{i=1}^{N} f_i(Z_{t_i-t_0})\right] + \varepsilon \prod_{i=1}^{N} \|f_i\|_\infty.$$

On the other hand,

$$\left|\mathbf{E}_{\phi m}\left[\prod_{i=1}^{N} f_i(Z_{t_i-t_0})\right] - \mathbf{E}_\mu\left[\prod_{i=1}^{N} f_i(Z_{t_i})\right]\right| < \varepsilon \prod_{i=1}^{N} \|f_i\|_\infty.$$

From these estimates we conclude that

$$\lim_{j \to \infty} \mathbf{E}^{k_j}_{\mu_{k_j}}\left[\prod_{i=1}^{N} f_i(Y^{k_j}_{t_i})\right] = \mathbf{E}_\mu\left[\prod_{i=1}^{N} f_i(Z_{t_i})\right].$$

However, we know that

$$\lim_{j \to \infty} \mathbf{E}^{k_j}_{\mu_{k_j}}\left[\prod_{i=1}^{N} f_i(Y^{k_j}_{t_i})\right] = \mathbf{E}\left[\prod_{i=1}^{N} f_i(Y_{t_i})\right].$$

This proves that $(Y, \mathbf{P})$ has the same finite dimensional distributions as those for the reflected Brownian motion $Z$ with initial distribution $\mu$. This completes the proof of the theorem. $\square$

REMARK 5.7. By [7], Proposition 3.10.4, Theorems 5.5 and 5.6 hold also for step-process approximation $\widehat{Y}^k$ defined as

$$\widehat{Y}_t^k := Y_{[2^k t]2^{-2k}}^k, \qquad t \geq 0,$$

but with the Skorokhod space $\mathbf{D}([0,1], \mathbb{R}^n)$ in place of the continuous function space $C([0,1], \mathbb{R}^n)$.

Now we turn to myopic processes $X^k$, defined at the beginning of Section 5.

LEMMA 5.8. *Suppose that either* (i) $\mu_k = m_k$ *for every* $k \geq 1$; *or* (ii) *there is a compact subset $K$ of $D$ such that $\mu_k$ is a measure on $D$ with* $\sup_{k \geq 1} \mu_k(D) < \infty$ *and* $\mu_k(D \setminus K) = 0$ *for all $k \geq 1$. Then the laws of* $\{X^k, \mathbf{P}_{\mu_k}^k, k \geq 1\}$ *are tight in the space $C([0,1], \mathbb{R}^n)$.*

PROOF. (i) Since, by Lemma 5.2, the laws of $\{Y^k, \mathbf{P}_{\mu_k}^k, k \geq 1\}$ are tight in the space $C([0,1], \mathbb{R}^n)$ and $Y_{j2^{-k}}^k = X_{j2^{-k}}^k$, we have, in particular, for every $\varepsilon > 0$ and $\delta > 0$, that there are $\rho > 0$ and $N \geq 1$ so that

$$(5.16) \qquad \mathbf{P}_{\mu_k}^k\left(\sup_{0 \leq i,j \leq 2^k : |i-j|2^{-k} \leq \rho} |X_{i2^{-k}}^k - X_{j2^{-k}}^k| > \varepsilon/3\right) < \delta/2$$

for every $k \geq N$.

Recall the oscillation operator $\mathrm{osc}_\rho[s,t]$ from the proof of Lemma 5.2 for the process $X^k$. As

$$\mathrm{osc}_\rho[0,1] \leq 2 \sup_{j \in \{1,\ldots,2^k\}} \mathrm{osc}_\rho[(j-1)2^{-k}, j2^{-k}]$$
$$+ \sup_{0 \leq i,j \leq 2^k : |i-j|2^{-k} \leq \rho} |X_{i2^{-k}}^k - X_{j2^{-k}}^k|,$$

we have

$$\mathbf{P}_{m_k}^k(\mathrm{osc}_\rho[0,1] > \varepsilon) \leq \mathbf{P}_{m_k}^k\left(\sup_{j \in \{1,\ldots,2^k\}} \mathrm{osc}_\rho[(j-1)2^{-k}, j2^{-k}] > \varepsilon/3\right)$$
$$(5.17)$$
$$+ \mathbf{P}_{m_k}^k\left(\sup_{0 \leq i,j \leq 2^k : |i-j|2^{-k} \leq \rho} |X_{i2^{-k}}^k - X_{j2^{-k}}^k| > \varepsilon/3\right).$$

Since $m_k$ is the reversible measure for the Markov chain $\{X^k_{j2^{-k}}, j = 0, 1, \ldots, 2^k\}$, for $k$ with $2^{-k} \leq \rho$,

$$
\begin{aligned}
\mathbf{P}^k_{m_k}&\left(\sup_{j \in \{1,\ldots,2^k\}} \mathrm{osc}_\rho[(j-1)2^{-k}, j2^{-k}] > \varepsilon/3\right) \\
&\leq \sum_{j=1}^{2^k} \mathbf{P}^k_{m_k}(\mathrm{osc}_\rho[(j-1)2^{-k}, j2^{-k}] > \varepsilon/3) \\
&= 2^k \mathbf{P}^k_{m_k}(\mathrm{osc}_\rho[0, 2^{-k}] > \varepsilon/3) \\
&\leq 2^k m(D) \mathbf{P}_0\left(\sup_{s,t \in [0,2^{-k}]} |X_s - X_t| > \varepsilon/3\right),
\end{aligned}
$$
(5.18)

which tends to zero as $k \to \infty$. Here $\mathbf{P}_0$ is the law of Brownian motion $X$ in $\mathbb{R}^n$ starting from the origin. Thus, (5.16)–(5.18) imply that there is $N_1 \geq N$ such that

$$\mathbf{P}^k_{m_k}(\mathrm{osc}_\rho[0, 1] > \varepsilon) < \delta \qquad \text{for every } k \geq N_1.$$

This proves that the laws of $\{X^k, \mathbf{P}^k_{m_k}, k \geq 1\}$ are tight in the space $C([0,1], \mathbb{R}^n)$.

(ii) can be established from the tightness of $\{X^k, \mathbf{P}^k_{m_k}, k \geq 1\}$ in the space $C([0,1], \mathbb{R}^n)$ by almost the same argument as that for the proof of Lemma 5.2(ii). So we omit the details here. □

Theorem 5.1 is a special case of the following with $\mu_k = \mu = \delta_{\{x_0\}}$, the Dirac measure concentrated at $x_0 \in D$.

THEOREM 5.9. *Suppose that either* (i) $\mu_k = m_k$ *for every* $k \geq 1$ *and* $\mu = m_D$; *or* (ii) $\{\mu_k, k \geq 1\}$ *is a sequence of measures on $D$ with* $\sup_{k \geq 1} \mu_k(D) < \infty$, $\mu_k(D \setminus K) = 0$ *for some compact subset $K$ of $D$ and all $k \geq 1$, and $\mu_k$ converge weakly on $D$ to a measure $\mu$. Then the laws of $\{X^k, \mathbf{P}^k_{\mu_k}, k \geq 1\}$ converge weakly in the space $C([0,1], \mathbb{R}^n)$ to the reflected Brownian motion on $D$ with the initial distribution $\mu$.*

PROOF. By Lemma 5.8, the laws of $\{X^k, \mathbf{P}^k_{\mu_k}, k \geq 1\}$ are tight in the space $C([0,1], \mathbb{R}^n)$. Let $(Z, \mathbf{P})$ be any of its subsequential weak limits. Clearly, $Z_0$ has distribution $\mu$. Let $Y$ be reflected Brownian motion on $D$ with initial distribution $\mu$. Denote by $\mathbb{Q}_2$ all the dyadic rational numbers. Since $X^k_t = Y^k_t$ for $t \in [0, 1]$ of the form $j2^{-k}$, we have from Theorems 5.5 and 5.6 that $\{Z_t, t \in \mathbb{Q}_2 \cap [0, 1], \mathbf{P}\}$ and $\{Y_t, t \in \mathbb{Q}_2 \cap [0, 1]\}$ have the same finite dimensional distributions. Since both $Z$ and $Y$ are continuous processes, $(Z, \mathbf{P})$ must have the same distribution as $Y$. That is, $(Z, \mathbf{P})$ is the reflected Brownian motion on $D$ with the initial distribution $\mu$. This proves the theorem. □

DEPARTMENT OF MATHEMATICS
UNIVERSITY OF WASHINGTON
SEATTLE, WASHINGTON 98195
USA
E-MAIL: burdzy@math.washington.edu
zchen@math.washington.edu